\documentclass{article}
\usepackage{blindtext}
\usepackage[a4paper, total={6in, 8in}]{geometry}

\usepackage{graphicx}
\usepackage{subfig}
\graphicspath{{./figures/}}
\usepackage{dcolumn}
\usepackage{bm}
\usepackage[mathlines]{lineno}
\usepackage[export]{adjustbox}
\usepackage{color}
\usepackage{tabularx,array,booktabs}
\usepackage{caption}
\usepackage{mathtools}
\usepackage{amsmath,amssymb}
\usepackage{soul}
\usepackage{varioref}
\usepackage{algorithm,algpseudocode}
\usepackage{ulem}
\usepackage{siunitx}
\usepackage{comment}
\usepackage{upgreek}
\usepackage{url}

\usepackage{calc}
\usepackage{multirow,multicol}

\newcommand{\crl}[1]{\textcolor{red}{#1}}

\newcommand{\bb}[1]{\mathbf{#1}}
\newcommand{\bs}[1]{\boldsymbol{#1}}
\newcommand{\R}{\mathbb R}
\newcommand{\C}{\mathbb C}
\newcommand{\LL}{\mathcal L}
\newcommand{\M}{\mathcal M}
\newcommand{\N}{\mathcal N}
\newcommand{\A}{\mathcal A}
\newcommand{\U}{\mathcal U}
\newcommand{\V}{\mathcal V}
\newcommand{\W}{\mathcal W}
\newcommand{\Y}{\mathcal Y}
\newcommand\norm[1]{\left\lVert#1\right\rVert}

\newcommand{\Matlab}{Matlab\textsuperscript{\textregistered}}

\DeclareMathOperator{\diag}{diag}

\begin{document}

\title{Optimal Strategies to Steer and Control Water Waves}

\date{}

  \author{Sebastiano Cominelli \thanks{Sebastiano Cominelli is a PhD Candidate at Politecnico di Milano, Department of Mechanical Engineering, Milano 20133, Italy (e-mail: sebastiano.cominelli@mail.polimi.it), Corresponding author.} \and
    Carlo Sinigaglia\thanks{Carlo Sinigaglia is a PhD Candidate at Politecnico di Milano, Department of Mechanical Engineering, Milano 20133, Italy (e-mail: carlo.sinigaglia@polimi.it)}  \and
  Davide E. Quadrelli \thanks{Davide E. Quadrelli is a PhD Candidate at Politecnico di Milano, Department of Mechanical Engineering, Milano 20133, Italy (e-mail: davidee.quadrelli@polimi.it).}\and
    Francesco Braghin\thanks{Prof. Francesco Braghin is Full Professor at Politecnico di Milano, Department of Mechanical Engineering, Milano 20133, Italy (e-mail: francesco.braghin@polimi.it). }
  }

\maketitle
\begin{abstract}
\noindent In this paper, we propose a novel approach for controlling surface water waves and their interaction with floating bodies. We consider a floating target rigid body surrounded by a control region where we design three control strategies of increasing complexity: an active strategy based on controlling the pressure at the air-water interface and two passive strategies where an additional controlled floating device is designed.
We model such device both as a membrane and as a thin plate and study the effect of this modelling choice on the performance of the overall controlled system. We frame this problem as an optimal control problem where the underlying state dynamics is represented by a system of coupled partial differential equations describing the interaction between the surface water waves and the floating target body in the frequency domain. An additional intermediate coupling is then added when considering the control floating device. The optimal control problem then aims at minimizing a cost functional which weights the unwanted motions of the floating body. A system of first-order necessary optimality conditions is derived and numerically solved using the finite element method. Numerical simulations then show the efficacy of this method in reducing hydrodynamic loads on floating objects.

\end{abstract}



\maketitle
\section{Introduction}
\label{sec:intro}

\noindent The idea of controlling water waves propagation through the design of suitable active or passive devices is mainly inspired by the cloaking theory originally developed by \cite{leonhardt2006optical,pendry2006controlling} for electromagnetic waves. The main theoretical tool adopted to render obstacles invisible to probes measuring the corresponding field consists of Transformation Theory (TT) whose output is a distribution of material properties which modify the wave propagation in the medium \cite{Kadic2016a}. Over the last two decades, similar techniques have been developed and adapted in a variety of different physical domains whose dynamical equations share a common mathematical structure, e.g.\ acoustics \cite{cummer2007one}, elasticity \cite{norris2011elastic}, heat conduction \cite{schittny2013experiments}, and water waves. 
In particular, Farhat et Al.\ \cite{farhat2008broadband} first applied TT to water waves, they designed a metamaterial for reducing the backscattering of a rigid obstacle irradiated by surface waves.
Later, \cite{newman2014cloaking}, followed by \cite{zhang2018quasi,zhang2020reduction}, studied the benefits given by many floating cylinders surrounding a fixed obstacle. \\
A different approach was proposed by \cite{porter2014cloaking,zareei2015cloaking}, the water waves are steered around a rigid target by choosing a proper shape of the surrounding seabed.
Then, \cite{dupont2016cloaking} 
proposed a structure made of many masses piercing the fluid that, coupled with a proper sea bed shape, can limit water actions on a submersed structure fixed to the seabed.

However, these methods cannot be adopted in case of floating devices, which are not rigidly connected to the ground. In \cite{alam2012broadband} was proposed to modify the seabed shape for shielding floating objects from gravity waves, but in practical cases, like wind turbines or other plants, the installation of those objects mostly occurs in deep sea zones where the sea bed shape is of little influence for wave propagation. \\
Inspired by modeling techniques for floating ice in oceans \cite{keller1953water,fox1990reflection}, a more realistic device was theorized in \cite{zareei2016cloaking}.
It consists of a thin floating plate that surrounds a circular infinite cylinder in constant depth environment and whose material properties are obtained by applying the conformal mapping method \cite{leonhardt2006optical}. However, their solution is based on assumptions that should be relaxed in case of floating objects since the obstacle is assumed to be fixed and extended until the seabed. Furthermore, the designed material properties require the plate to be anisotropic and inhomogeneous, which may be difficult to manufacture in practice.
Recently, \cite{loukogeorgaki2019minimization} proposed to reduce the oscillations of a floating cylinder by optimizing the thickness of a homogeneous annular plate that floats around the cylinder; differently, \cite{iida2023water} considered a composite plate made of many concentric homogeneous layers. Both strategies are based on the analytical solution for axisymmetric floating objects.

In this paper, we formulate an Optimal Control Problem (OCP) which allows to take into account hydrodynamic interactions between water and target floating body. The OCP formulation is also able to tackle complex geometries and practical constraints on the material properties. Furthermore, it is relatively straightforward to encode different objectives in the OCP formulation. In this way, we are able to design devices close to real-world applications.
We choose control mechanisms that act on a region surrounding the obstacle. In particular, we investigate both active and passive control strategies that interact with the water surface.\\
The rigid-body motion of the floating structure can be measured by a cost function, so the problem is addressed as an Optimal Control Problem (OCP) constrained by a system of coupled Partial Differential Equations (PDEs) that govern the water and the control device dynamics.\\
In search of an optimal solution, we derive first order necessary conditions by using a Lagranian approach, see e.g. \cite{manzonioptimal}, which are then solved by an iterative, gradient-based optimization algorithm for some relevant test cases.\\
A remarkable advantage of this approach is that problems of arbitrary geometries can be considered, which is of great importance in engineering applications. Moreover, similar to \cite{cominelli2022design}, where a method for designing acoustic cloaks by solving PDE-constrained OCPs is proposed, we expect
that a narrowband high-performance device will be achieved, and this is a good premise for applications involving swells, which are steady-state full developed water waves generated by distant storms, and are characterized by narrowband spectra \cite{hasselmann1973measurements}.

In the following we briefly review the paper organization.
In Section~\ref{sec:physics} we describe the physical model adopted together with its main assumptions, then in Section~\ref{sec:pressure controlled} we propose an active control strategy that modifies the pressure on the air/water interface surrounding the floating target.

The resulting OCP is linear-quadratic and can be solved efficiently as a large linear system, additionally, fast solving techniques are available for real-time applications, exploiting, for example, Model Order Reduction strategies \cite{sinigaglia2022fast}.
This control idea is well suited for floating systems such as Floating Production Storage and Offloading (FPSO), where energy can be spent for active control, conversely this is hardly the case of floating turbines where a passive mechanism is sought instead.
After treating the active linear problem, we consider passive strategies in Sections~\ref{sec:membrane controlled} and \ref{sec:plate controlled}. In particular, we analyse two control mechanisms that consist of covering a portion of the water surface around the turbine with a floating elastic membrane or a plate, respectively. The two models describing the coupled dynamics have been originally proposed for studying the effect of floating ice on water waves \cite{keller1953water,fox1990reflection}, and they now come in handy as constraints for the OCP. \\
In this case, the control is tuned by modifying stiffness and inertia properties of the floating device: by increasing the local stiffness, water waves experiment a speed up, conversely they are slowed down as the inertia increases. Note that, differently from \cite{zareei2016cloaking}, this framework guarantees the isotropy of both the membrane and the plate, as a consequence , we obtain a simpler device to manufacture.




\section{Problem Statement}
\label{sec:physics}
\begin{figure}
    \centering
    \subfloat[][]
    {\includegraphics[width=.35\textwidth]{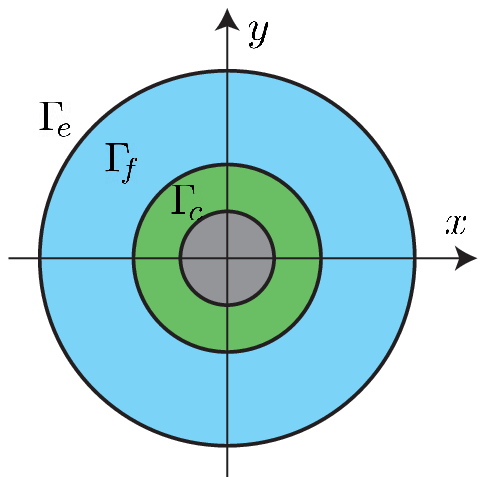}}
    \quad
    \subfloat[][]
    {\includegraphics[width=.55\textwidth]{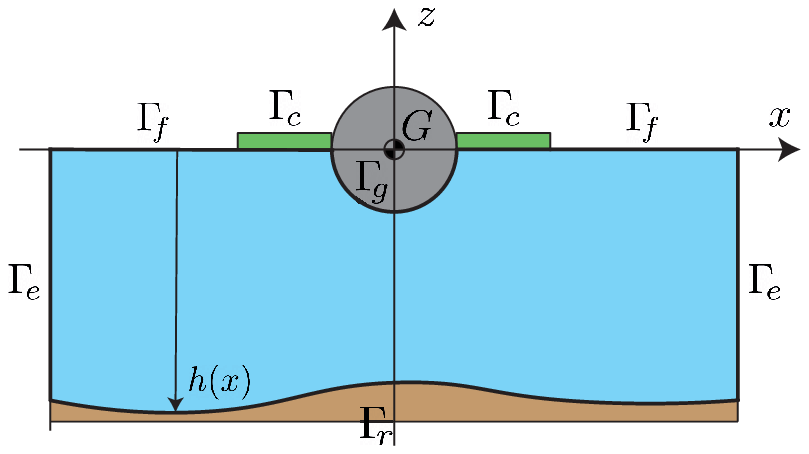}}
    \caption{top (a) and side (b) views of the computational domain. The floating body is a sphere with mass density half than water, the green region represents the control surface.
    }
    \label{fig:computational domain}
\end{figure}

Throughout this section, we briefly sum up the main modelling equations applied for the rest of the paper giving a summary of their derivation; this turns out to be useful for Sections~\ref{sec:pressure controlled}, \ref{sec:membrane controlled} and \ref{sec:plate controlled} in particular, where the different control actions considered lead to modifications of the system dynamics.\\
For a detailed derivation, the reader is referred to the monograph \cite{mei2005theory}. The main assumptions we make are that the flow is inviscid, irrotational, and characterized by small wave amplitude.

Let us consider the volume of water $\Omega\subset\R^3$ depicted in Figure~\ref{fig:computational domain} representing top and side views of a cylindrical portion of an ocean environment. Its disjoint boundaries are $\Gamma_f$, $\Gamma_c$, $\Gamma_r$, $\Gamma_g$ and $\Gamma_e$. $\Gamma_f$ and $\Gamma_c$ are the free and the controlled parts of air-water interface respectively, described by the time-varying surface $z=\zeta(x,y,t)$ and with equilibrium position $z=0$; $\Gamma_r$ is the sea bottom surface described by the function $h=h(x,y)$. $\Gamma_g$ is the wetted surface of the target floating body and its position depends on the body motion; $\Gamma_e$ is the artificial boundary introduced for computational reasons.
    
By assuming the variation in water density insignificant over the temporal and spatial scales of interest, the continuity and Navier-Stokes equations adequately describing the fundamental conservation laws of mass and momentum are respectively:
\begin{align}\label{eq:mass conservation}
    \nabla\cdot\bf u &= 0,
    \\ \label{eq:momentum conservation}
    \bb u_t +\bb u\cdot\nabla\bb u &= -\nabla\left(\frac{P}{\rho}+gz\right)+\upsilon\Delta\bb u,
\end{align}
inside the computational domain $\Omega$, where $\bb u=\bb u(\bb x,t) \in \mathbb{R}^{3}$ is the velocity vector field and $\bb u_t$ its partial time derivative while $\bb{u} \cdot \nabla \bb u$ denotes the convective terms using the classical fluid dynamics notation, i.e. $\displaystyle (\bb{u} \cdot \nabla \bb u)_{i} = \sum_{j=1}^{3} u_j \frac{\partial u_i}{\partial x_j}$. Furthermore, $P=P(\bb x,t)$ is the pressure field, $\rho=\SI{1030}{\kilo\gram\per\cubic\meter}$ the water density, $g=\SI{9.81}{\meter\per\square\second}$ the gravitational acceleration, $\upsilon$ the kinematic viscosity coefficient; the vector $\bb x = (x,y,z)$ is referred with the $z$ axis pointing vertically upward.
Under the assumption of inviscid ($\upsilon = \SI{0}{\square\meter\per\second}$) irrotational flow, the velocity field can be expressed as the gradient of a scalar potential $\Phi(\bb x,t)$, see e.g. \cite{mei2005theory} 
\begin{equation}\label{eq:irrotational flow}
    \bb u = \nabla\Phi.
\end{equation}
So, by merging \eqref{eq:mass conservation} and \eqref{eq:irrotational flow}, $\Phi$ must satisfy the Laplace equation
\begin{equation}\label{eq:Poisson}
    \Delta\Phi = 0.
\end{equation}
Under the same assumptions, the well known unsteady Bernoulli equation can be proved starting from  \eqref{eq:mass conservation} and \eqref{eq:momentum conservation} we obtain
\begin{align}\label{eq:Bernoulli}
    -\frac{P}{\rho} &= gz +\Phi_t +\frac{1}{2}\left\|\nabla\Phi\right\|^2 &\text{in }\Omega.
\end{align}
Finally, if the ratio $\varepsilon\coloneqq2\pi \frac{A}{\lambda}$ between the wave amplitude $A$ and the wavelength $\lambda$ is much smaller than one, that is $\varepsilon \ll 1$, the quadratic term in \eqref{eq:Bernoulli} can be neglected obtaining the linear relation
\begin{align}\label{eq:lin Bernoulli}
    -\frac{P}{\rho} &= gz +\Phi_t  &\text{in }\Omega.
\end{align}

Throughout this paper, we develop control strategies that work indefinitely in time for systems characterized by specific frequencies, hence the frequency domain approach is adopted for reducing the computational effort in view of numerical solutions. Adopting the usual notation, let us introduce the following variables:
\begin{align}
    \begin{pmatrix}
            \Phi(\bb x,t) \\ \zeta(x,y,t) \\ P(\bb x,t)
    \end{pmatrix}
    & = \Re\left\{
    \begin{pmatrix}
            \phi(\bb x) \\ \eta(x,y)  \\ p(\bb x)
    \end{pmatrix}
    e^{j\omega t}
    \right\}.
\end{align}
where $\phi$, $\eta$ and $p$ are complex valued functions and $\omega$ is the angular frequency.
Thus, equations~\eqref{eq:Poisson} and \eqref{eq:lin Bernoulli} become, respectively
\begin{align}\label{eq:laplace freq}
    \Delta\phi &= 0,
    \\ \nonumber
    -\frac{p}{\rho} &= gz+j\omega\phi.
\end{align}

\subsection{Boundary conditions}
We now briefly describe the remaining conditions that define a boundary value problem for the velocity potential $\phi$.\\
The seabed is considered as an infinitely rigid boundary; then, from equation~\eqref{eq:irrotational flow}, we can state that the velocity component normal to $\Gamma_r$ is null, i.e.\
\begin{align}\label{eq:BC on G_r}
    \phi_n &= 0 &\text{on }\Gamma_r,
\end{align}
where the subscript $n$ stands for the derivative along the outgoing normal of the boundary.

Let us now consider the boundaries $\Gamma_f$ and $\Gamma_c$: they belong to the surface ${z=\zeta(x,y,t)}$, whose shape depends both on space and time. However, under the assumption of small wave amplitude ($\varepsilon\ll1$), it can be approximated up to the first order as the plane $z=0$.\\
Using a first order expansion of the surface $z - \eta=0$ and considering that $\varepsilon\ll1$, one can obtain the following kinematic relationship
\begin{align} \label{eq:kinematic}
    \phi_n &= j\omega\eta &\text{on }\Gamma_f\cup\Gamma_c
\end{align}
that describes the continuity of velocity while neglecting the convective terms, see \cite{mei2005theory} for a detailed derivation. In addition, on the air-water interface, the Bernoulli equation \eqref{eq:lin Bernoulli} states that
\begin{align}\label{eq:Bernoulli G_f}
    -\frac{p}{\rho} &= g\eta + j\omega\phi  &\text{on }\Gamma_f\cup\Gamma_c.
\end{align}
So, on the free surface $\Gamma_f$, the boundary condition on $\phi$ can be obtained by merging the last two expressions. Multiplying \eqref{eq:Bernoulli G_f} by $j\omega$ and using \eqref{eq:kinematic}, we obtain a boundary condition for $\phi$ on $\Gamma_f$:
\begin{align}\label{eq:pressure on G_f}
    g\phi_n -\omega^2\phi  &= -j\omega\frac{p}{\rho} & \text{on }\Gamma_f,
\end{align}
where the right-hand-side represents a forcing term coming from the environment, e.g.\ from wind.
Conversely, the dynamics holding on the control surface $\Gamma_c$ depends on the kind of control adopted. In this case, the kinematic relation \eqref{eq:kinematic} remains valid, while \eqref{eq:Bernoulli G_f} changes according to the cases analyzed in Sections~\ref{sec:pressure controlled}, \ref{sec:membrane controlled} and \ref{sec:plate controlled}. We define an operator $E$ which encodes the dynamic equilibrium that fluid velocity potential $\phi$, vertical displacement $\eta$ and control action $u$ shall satisfy:
\begin{align*}
    E(\phi,\eta,u) &= 0 &\text{on }\Gamma_c.
\end{align*}

For what concerns the wet surface of the floating body, $\Gamma_g$, two conditions have to be imposed: kinematic constraints on $\Gamma_g$ and the balance of forces on the body. Even though the boundary moves following the motion of the body, up to the first order approximation it can be assumed fixed in its equilibrium position for solving in a simplified way the boundary value problem \cite{mei2005theory}, such that the following two equations describing the linearized motion can be derived:
\begin{align}\label{eq:body kinematics}
    &\phi_n = j\omega\{\bb n\}^\top \{\bb X\}&\text{on }\Gamma_g
    \\ \label{eq:body dynamics}
    &[K-\omega^2M]\{\bb X\} = -j\omega\rho\int_{\Gamma_g}\phi\{\bb n\}d\Gamma + \{\bb f\}
\end{align}
where $\{\bb X\} = (\bb x^b, \bs\theta^b)\in\C^6$ is the column vector collecting the six degrees of freedom of the floating body in frequency domain: the three displacements $\bb x^b\in\C^3$ and the three rotation angles $\bs\theta^b\in\C^3$ of the body with respect to a fix point $G$; ${\{\bb n\} = (\bb n, (\bb x- \bb x_G)\times \bb n)}$, ${\{\bb n\}\colon \R^3\to \R^6}$ is the generalized normal to the surface, with $\bb n = \bb n(\bb x)$ the outgoing normal and $\bb x_G$ the spatial coordinates of $G$. We indicate with $\times$ the cross product between two vectors. \\
The vector $\{\bb f\}\in\C^6$ collects the external loads acting on the body and allows to consider the effect of forces due to wind and catenary mooring lines \cite{mei2005theory}; for the sake of simplicity, we suppose $\{\bb f\} = \bb0$ for the rest of this paper. The floating body dynamics \eqref{eq:body dynamics} depends on the stiffness and mass matrices $K$, $M\in\R^{6\times6}$ respectively, they are defined as:
\begin{equation}\label{eq:M and K matrices}
\begin{aligned}
    K =&
    \rho g\begin{bmatrix}
    0  &  0  &  0     &    0               &    0              &  0             \\
    0  &  0  &  0     &    0               &    0              &  0             \\
    0  &  0  &  \A     &  I^\A_2           & -I^\A_1           &  0             \\
    0  &  0  &  I^\A_2 &  I_{22}^\A+I_3^V  & -I^\A_{12}        &  0             \\
    0  &  0  & -I^\A_1 & -I^\A_{21}        &  I_{11}^\A+I_3^V  &  0             \\
    0  &  0  &  0     &    0               &    0              &  0
    \end{bmatrix} ,
    \\
    M = &
    \begin{bmatrix}
    M^b  &  0  &  0  &          0          &          0          &  0             \\
    0  &  M^b  &  0  &          0          &          0          &  0             \\
    0  &  0  &  M^b  &          0          &          0          &  0             \\
    0  &  0  &  0  &  I^b_{22}+I^b_{33}  &      -I^b_{21}      &  -I^b_{31}     \\
    0  &  0  &  0  &  -I^b_{12}          &  I^b_{33}+I^b_{11}  &  -I^b_{32}     \\
    0  &  0  &  0  &  -I^b_{13}          &      -I^b_{32}      &  I^b_{11}+I^b_{22} 
    \end{bmatrix} ,
\end{aligned}
\end{equation}
where $\A$ is the area of $S^\A$, with $S^\A$ the cross-section of the body with respect to the plane $z=0$, $M^b$ the body mass, $I^\A$ and $I^b$ are the first and second moments of inertia with respect of the surface $S^\A$ and the body volume $V^b$ respectively, and $I^V$ is the moments of inertia of the submersed volume $V$, i.e.:
\begin{align*}
    I_i^\A &= \int_{S^\A}{\big( \bb x - \bb{x}_G\big)_i \,dS}
    \\
    I_{ij}^\A &= \int_{S^\A}{ \big(\bb x - \bb x_G\big)_i\big(\bb x - \bb x_G\big)_j\,dS}
    \\
    I_i^V &= \int_{V}{\big( \bb x - \bb{x}_G\big)_i \,dV}
    \\
    I_{ij}^b &= \int_{V^b}{ \big(\bb x - \bb x_G\big)_i\big(\bb x - \bb x_G\big)_j \,dm}
\end{align*}

For computational reasons, the domain has been truncated generating a fictitious cylindrical boundary $\Gamma_e$ on which an absorbing condition must be imposed for avoiding artificial reflections. For the sake of simplicity, we consider a first order radiation condition (see e.g.\ \cite{bai1972variational,hsu2003second}), which is able to absorb waves with a small incidence angle on $\Gamma_e$ with respect to the normal:
\begin{align}\label{eq:radiating sc}
    &\phi^s_n + \alpha\phi^s = 0 &  \text{on } \Gamma_e
\end{align}
where $\alpha=jk+\frac{1}{2R}$; $\phi^s$ is the scattered potential field with respect to the incident one $\phi^i$, characterized by a single angular frequency $\omega$ and wavenumber $k$; $R$ is the base radius of the cylinder $\Gamma_e$. Note that since the floating obstacle is in the middle of the computational domain and the absorbing boundary is a cylinder surrounding the body, scattering is expected to come nearly orthogonal to $\Gamma_e$, thus leading to an acceptable numerical approximation for a computational domain sufficiently large relatively to the floating body and the control mechanism.\\
The incident field $\phi^i$ is the analytical solution of a wave propagating in a domain without obstacles and whose depth is constant, i.e.\ $h(x,y) \equiv h_0$:
\begin{equation}\label{eq:ana sol}
    \phi^i = j\frac{gA}{\omega}\frac{\cosh{k(z+h_0)}}{\cosh{kh_0}}e^{j\bb k\cdot \bb x},
\end{equation}
where $\bb k$ is the wave vector and $A$ is the wave amplitude. Wave number $k=\lvert\bb k\lvert$ and circular frequency $\omega$ must satisfy the dispersion relation
\begin{equation}\label{eq:dispersion}
    \omega^2 = g k\tanh{kh_0}.
\end{equation}
In other words, \eqref{eq:ana sol} is the analytical solution to the potential equation~\eqref{eq:laplace freq}, the free surface equilibrium~\eqref{eq:pressure on G_f} with $p=0$ and the sea bed condition~\eqref{eq:BC on G_r} in case there are no floating obstacles and the sea depth is constant and equal to $h_0$. Again, the reader is referred to, e.g., \cite{mei2005theory} for a detailed derivation. \\
The radiation condition~\eqref{eq:radiating sc} holds for the scattered field only, so we have to reformulate all the above equations in terms of $\phi^s$ making use of $\phi = \phi^s + \phi^i$, since $\phi$ is considered as the total velocity potential.
We obtain
the following linear elliptic PDE coupled with the rigid body dynamics and the pressure equilibrium on $\Gamma_c$:
\begin{align}
&\begin{dcases}
\label{eq:potential dynamics}
    -\Delta\phi^s = 0   & \text{in }\Omega
    \\
    \phi^s_n = 0  &\text{on }\Gamma_{r}
    \\
    \phi^s_n -\frac{\omega^2}{g}\phi^s = 0&\text{on }\Gamma_f
    \\
    \phi^s_n = j\omega\eta - \phi^i_n &\text{on }\Gamma_{c}
    \\
    \phi^s_n = j\omega\{\mathbf n\}^\top \{\mathbf X\} - \phi^i_n&\text{on } \Gamma_g
    \\
    \phi^s_n + \alpha\phi^s = 0& \text{on }\Gamma_e
\end{dcases}
\\[5pt] \label{eq:body dynamics2}
&\quad\big(K-\omega^2M\big)\{\bb X\} = -j\omega\rho\int_{\Gamma_g}(\phi^s+\phi^i)\{\bb n\} \,d\Gamma
\\[5pt] \nonumber
&\quad E(\phi,\eta,u) = 0 \qquad\qquad\qquad\text{on }\Gamma_c
\end{align}
Note that, for the sake of simplicity, the wind pressure on $\Gamma_f$ is assumed to be null so the incoming wave is considered as generated away from the region of interest, that is a consistent hypothesis in case of swells. Also, the sea bed is considered of constant depth.
In this setting, the scattered field $\phi^s$ includes the perturbation with respect to $\phi^i$ given by the floating body and due to any control action we apply.

\section{Active pressure control}
\label{sec:pressure controlled}

In this Section we formulate and solve the active control problem where the surface pressure around the obstacle is considered as control mechanism. The control space is denoted by $\mathcal{U}$ and it is selected as $L^{\infty}(\Gamma_c,\C)\cap H^1(\Gamma_c,\C)$. We define the cost functional $J\colon\C^6\times\U\to\R$ as
\begin{align}\label{eq:cost}
    J\big(\{\bb X\},u\big) \coloneqq \frac{1}{2}\{\bb X\}^\dag C\{\bb X\} + \frac{1}{2}\norm{u}_{\U}^2,
\end{align}
where $\{\bb X\}^\dag=\{\bar{\bb X}\}^\top$ stands for the hermitian of $\{\bb X\}$ and $C\in\R^{6\times6}$ is a positive definite weighting matrix; the norm $\norm{u}_\U$ on the space $\U$ is defined as the weighted norm
\begin{equation*}
    \norm{u}^2_\U\coloneqq
    \left(\bar u, u\right)_\U=
    \alpha_u\norm{u}^2_{L^2(\Gamma_c)} + \beta_u\norm{\nabla u}^2_{L^2(\Gamma_c)},
\end{equation*}
where $\alpha_u,\beta_u>0$ are regularization parameters that allow us to limit separately the control effort and the spatial control gradient respectively. In particular, $\alpha_u$ limits the control effort, while $\beta_u$ imposes a soft constraint on the control gradient. \\
Note that the cost functional involved in the OCP is the same independently on the control mechanism since through all the paper we aim at controlling the water flow such that the floating body motion is minimized.

The first scenario here considered is that of an active source applying a pressure field on the boundary $\Gamma_c$, surrounding the floating body. An exact solution can be obtained in this case without any iterative algorithm, due to the linear-quadratic nature of the problem. This fact makes the following discussion appealing for real-time applications, e.g.\ for Linear Quadratic Regulators \cite{manzonioptimal}.

Hence, we look for a surface pressure $p\lvert_{\Gamma_c} = u$, $u\colon\R^2\to\C$ such that the cost functional $J$ is minimized.
Since the only action on the surface $\Gamma_c$ is a pressure, the condition $E(\phi,\eta,u)=0$ can be obtained by the Bernoulli equation \eqref{eq:lin Bernoulli} evaluated on $\Gamma_c$:
\begin{align*}
    E(\phi,\eta,u) = g\eta + j\omega\phi + \frac{u}{\rho}.
\end{align*}
It can be merged with the kinematic relation~\eqref{eq:kinematic}
\begin{align*}
    \phi_n -\frac{\omega^2}{g}\phi  &= -j\frac{\omega}{\rho g}u & \text{on }\Gamma_c
\end{align*}
and
\begin{align*}
    \phi^s_n -\frac{\omega^2}{g}\phi^s  &= -j\frac{\omega}{\rho g}u & \text{on }\Gamma_c
\end{align*}
since $\phi_n^i-\frac{\omega^2}{g}\phi^i=0$ on $\Gamma_c$. Then the OCP reads
\begin{equation}\label{eq:OCP 1}
\begin{gathered}
    \min_{u\in\mathcal{U}_{ad}}{\tilde J} = J\Big(\{\bb X\}(u), u\Big)
        \\
    \text{s.t.} \quad
        \begin{aligned}
            &\qquad\begin{dcases}
                -\Delta\phi^s = 0   & \text{in }\Omega
                \\
                \phi^s_n = 0  &\text{on }\Gamma_{r}
                \\
                \phi^s_n -\frac{\omega^2}{g}\phi^s = 0&\text{on }\Gamma_f
                \\
                \phi^s_n -\frac{\omega^2}{g}\phi^s = -j\frac{\omega}{\rho g}u &\text{on }\Gamma_{c}
                \\
                \phi^s_n = j\omega\{\mathbf n\}^\top \{\mathbf X\} -\phi^i_n &\text{on } \Gamma_g
                \\
                \phi^s_n + \alpha\phi^s = 0 & \text{on }\Gamma_e
            \end{dcases}
            \\[5pt]
            &\qquad\quad  \big(K-\omega^2M\big)\{\bb X\} = -j\omega\rho\int_{\Gamma_g}(\phi^s+\phi^i)\{\bb n\}d\Gamma
        \end{aligned}
    \end{gathered}
\end{equation}
For the sake of simplicity, we choose $\U_{ad} = \U$, the existence and uniqueness of a solution is guaranteed
since $\alpha_u$, $\beta_u$ are strictly positive \cite{manzonioptimal}. We follow a Lagrange multiplier method to derive a system of first-order necessary conditions; see, e.g. \cite{manzonioptimal}. We define the Lagrangian $\LL\colon\W\times\W'\times\C^6\times\C^6\times\U\to\R$ as
\begin{multline}\label{eq:weak form pressure}
    \LL=
    J
    +\Re\Big\{       \int_\Omega{\nabla\bar\lambda\cdot\nabla\phi^s}
    -\frac{\omega^2}{g}\int_{\Gamma_f}{\bar\lambda\phi^s}
    -\frac{\omega^2}{g}\int_{\Gamma_c}{\bar\lambda\phi^s}
    +j\frac{\omega}{\rho g}\int_{\Gamma_c}{\bar\lambda u}
    -j\omega\int_{\Gamma_g}{\bar\lambda\{\bb n\}^\dag}\,\{\bb X\}
    \\
    +\int_{\Gamma_g}{\bar\lambda\phi^i_n}
    + \alpha\int_{\Gamma_e}{\bar\lambda\phi^s}
    +\{\bb Y\}^\dag 
    \Big[\big(K-\omega^2M\big)\{\bb X\} +j\omega\rho\int_{\Gamma_g}(\phi^s+\phi^i)\{\bb n\}\Big]
    \Big\}.
\end{multline}
where $\lambda\in\W'$ is the adjoint of the state $\phi^s$ and $\W'=H^1(\Omega,\C)$ its functional space, $\{\bb Y\}\in\C^6$ the adjoint of $\{\bb X\}$. With a slight abuse of notation, we omit the differentials inside integrals. Note that the real part of the scalar products is taken, but similar results can be obtained by considering the imaginary part instead.\\
By taking the Gate\^aux derivative of $\LL$ with respect to the state $\phi^s$, the derivative with respect to $\{\bb X\}$ , and by requiring they are null, we can obtain the optimality conditions on the adjoints:
\begin{gather*}
    \LL'_{\phi^s}[\varphi] = \frac{1}{2}    \Big\{      
    \int_\Omega{\nabla\bar\lambda\cdot\nabla\varphi}
    -\frac{\omega^2}{g}\int_{\Gamma_f\cup\Gamma_c}{\bar\lambda\varphi}
    + \alpha\int_{\Gamma_e}{\bar\lambda\varphi}
    + j\omega\rho\{\bb Y\}^\dag\int_{\Gamma_g}{\varphi\{\bb n\}}
    \Big\} = 0 \qquad\forall\varphi\in\W,
    \\
    \LL'_{\{\bb X\}} =     \frac{1}{2}\{\bb X\}^\dag C
    +\frac{1}{2}\Big\{
    -j\omega\int_{\Gamma_g}{\bar\lambda\{\bb n\}^\dag}
    +\{\bb Y\}^\dag \big(K-\omega^2M\big)
    \Big\} = \bb 0,
\end{gather*}
where Wirtinger's calculus rules have been applied for handling complex derivatives \cite{wirtinger1927formalen}.
The respective strong formulations are:
\begin{equation}\label{eq:pressure adj}
\begin{gathered}
    \begin{dcases}
    -\Delta\lambda = 0 &\text{in }\omega
    \\
    \lambda_n = 0 &\text{on }\Gamma_r
    \\
    \lambda_n - \frac{\omega^2}{g}\lambda= 0 &\text{on }\Gamma_f\cup\Gamma_c
    \\
    \lambda_n = j\omega\rho \{\bb n\}^\top\{\bb Y\} &\text{on }\Gamma_g
    \\
    \lambda_n +\bar\alpha\lambda = 0 &\text{on }\Gamma_e
    \end{dcases}
    \\[5pt]
    (K-\omega^2M)^\top\{\bb Y\} = -j\omega\int_{\Gamma_g}{\lambda\{\bb n\}\,d\Gamma} - C^\top\{\bb X\}
\end{gathered}
\end{equation}
Finally, the reduced cost gradient is obtained by taking the derivative of \eqref{eq:weak form pressure} with respect to $u$:
\begin{align}
    \LL'_u[\psi] &= \int_{\Gamma_c} \left( \frac{\alpha_u}{2}\bar u\psi +\frac{\beta_u}{2}\nabla\bar u\cdot\nabla\psi + j\frac{\omega}{2\rho g}\bar\lambda\psi\right)\,d\Gamma =0
    &\forall\psi\in\U_{ad};
\end{align}
note that the equality holds since $\U_{ad}=\U$. 

In search of a numerical solution, we rely on a Galerkin $\mathbb P_2$ Finite Element Method (FEM) for discretizing the optimization problem. Thus a triangulation $\mathcal T_h$ of $\Omega$ with characteristic size $h>0$ is considered and the discrete spaces $\W_h$ and $\W'_h$ of $\W$ and $\W'$ are defined respectively; accordingly, $\U_h$ is the discrete control space defined on the discrete counterpart of $\Gamma_c$. As a consequence, state, adjoint and control are discretized as
\begin{align*}
    \phi^s(\bb x) \approx \bs\varphi(\bb x)^\top \bs\phi        ,\qquad
    \lambda(\bb x) \approx \bs\varphi(\bb x)^\top \bs\lambda  ,\qquad 
    u(\bb x) \approx \bs\psi(\bb x)^\top \bb u,
\end{align*}
The discrete cost functional and the optimality necessary conditions are
\begin{align*}
    \text{Cost functional:} & \qquad
        J = \frac{1}{2}\{\bb X\}^\dag C\{\bb X\} + \frac{\alpha_u}{2}\bb u^\dag E_c \bb u + \frac{\beta_u}{2}\bb u^\dag A_c \bb u
    \\[5pt]
    \text{State:} & \qquad
    \begin{aligned}
        &\left[ A-\frac{\omega^2}{g}(C_f+C_c)+\alpha C_e \right] \bs\phi = j\omega K_g^\top\{\bb X\} + \bb f_g - j\frac{\omega}{\rho g} D_c\bb u
        \\
        &(K-\omega^2M)\{\bb X\}=-j\omega\rho (K_g\bs\phi +\bb g)
    \end{aligned}
    \\[5pt]
    \text{Adjoint state:} & \qquad
    \begin{aligned}
        &\left[ A - \frac{\omega^2}{g}(C_f+C_c) +\bar\alpha C_e \right]\bs\lambda = j\omega\rho K_g^\top\{\bb Y\} 
        \\
        &(K-\omega^2M)\{\bb Y\}=-j\omega K_g\bs\lambda - C^\top\{\bb X\}
    \end{aligned}
    \\[5pt]
    \text{Control gradient:} &\qquad
        \nabla_{\bb u} J = \left(\frac{\alpha_u}{2}E_c
        +\frac{\beta_u}{2}A_c\right)\bar{\bb u}
        + j\frac{\omega}{2\rho g} D_c^\top\bar{\bs\lambda} = \bb 0
\end{align*}
where $A\in\R^{n\times n}$ is the stiffness matrix on the domain $\Omega$; $C_f$, $C_c$ and $C_e\in\R^{n\times n}$ are the mass matrices on the boundaries $\Gamma_f$, $\Gamma_c$ and $\Gamma_e$ respectively. $D_c\in\R^{n\times l}$ is the control matrix and $A_c,E_c\in\R^{l\times l}$ are the stiffness and the mass matrices in $\Gamma_c$. The matrix $K_g\in\R^{6\times n}$ is defined as
\begin{equation*}
    (K_g)_{ij} \coloneqq \int_{\Gamma_g}\{\bb n\}_i\,\varphi_j\,d\Gamma;
\end{equation*}
$\bb f_g\in\C^n$ and $\bb g\in\C^6$ are forcing terms given by the background field, i.e.\
\begin{align*}
    \bb f_g = \int_{\Gamma_g} -\phi^i_n\,\bs\varphi \,d\Gamma
    &&
    \bb g = \int_{\Gamma_g} \phi^i\{\bb n\}\,d\Gamma
\end{align*}
The OCP \eqref{eq:OCP 1} is a linear quadratic one, as a consequence, it is also convex given that the weighting matrices in the cost functional are positive definite \cite{manzonioptimal} and the necessary conditions are also sufficient; then its global minimum can be obtained by solving the following linear system
\begin{align}\label{eq:pressure lin system}
        \begin{bmatrix}
            A_{tot}         & -j\omega K_g^\top & \bb 0         & \bb 0                 & j\frac{\omega}{\rho g} D_c \\
            j\omega\rho K_g & (K-\omega^2M)     & \bb 0         & \bb 0                 & \bb 0 \\
            \bb 0           & \bb 0            & A_{tot}^\dag  & -j\omega\rho K_g^\top & \bb 0 \\
            \bb 0           & C^\top             & j\omega K_g   & (K-\omega^2M)         & \bb 0 \\
            \bb 0           & \bb 0             & -j\frac{\omega}{2\rho g}D_c^\top &\bb0         & \frac{\alpha_u}{2}E_c + \frac{\beta_u}{2}A_c
        \end{bmatrix}
        \begin{pmatrix}
             \bs \phi \\
             \{\bb X\} \\
             \bs\lambda \\
             \{\bb Y\} \\
             \bb u
        \end{pmatrix}
        =
        \begin{pmatrix}
             \bb f_g \\ j\omega\rho\,\bb g \\ \bb 0 \\ \bb 0 \\ \bb 0
        \end{pmatrix},
\end{align}
where $A_{tot} = A-\frac{\omega^2}{g}(C_f+C_c)+\alpha C_e$. 

As a case study, we consider a floating sphere with homogeneous density $\rho_b = 0.5\rho$ such that, in the static equilibrium condition, the body is half submersed as shown in figure~\ref{fig:computational domain}. The center of mass is assumed to be in the origin, i.e.\ $\bb x_G = \bb 0$. Due to the spherical symmetry, the only non-null term of the $K$ matrix is $(K)_{33}$. This in turns implies that the interaction with water cannot result in a rotation of the body, since $(\bb x-\bb x_G)\times\bb n=\bb0$ on $\Gamma_g$ and water viscosity is considered null. Since the dynamics is linear, we assume a unitary wave amplitude of the incident wave \eqref{eq:ana sol} for the sake of simplicity.
All geometrical sizes are normalized with respect to the sphere diameter, so we choose a flat seabed \SI{2.5}{\meter} deep, i.e.\ $h(\bb x) = \SI{2.5}{\meter}$, and we limit the computational domain to a cylinder of radius \SI{4}{\meter}. The control surface is an annular region surrounding the sphere whose external radius is \SI{1}{\meter} long.

We choose the excitation wave period of $T=\SI{1.2}{\second}$ such that the normalized wavelength $\lambda\approx\SI{2.25}{\meter}$
results in a relevant excitation of the floating body. Note that, since the aim is not to cloak the floating obstacle in the usual sense, the wavelengths of interest are very different from the ones considered in usual cloaking problems; in particular, a very short $\lambda$ results in a small excitation of the target since the mean pressure on $\Gamma_g$ is negligible. Conversely, a wave characterized by a relatively long wavelength push on the floating body coherently.
The wave is assumed to propagate along the y-axis, i.e.\ $\bb k= (0, k, 0 )$.
Finally, we choose the cost $C$ to be the diagonal matrix
$C = \diag(1,1,1,H^2,H^2,H^2)$, where $H=\SI{1}{\meter}$ is the characteristic height of the floating object. Finally, the regularization parameters are chosen as $\alpha_u=\beta_u=\SI{1e-10}{}$.

\begin{figure}
    \centering
    \subfloat[][]
    {   \includegraphics[trim={40 0 25 0},width=.45\textwidth]{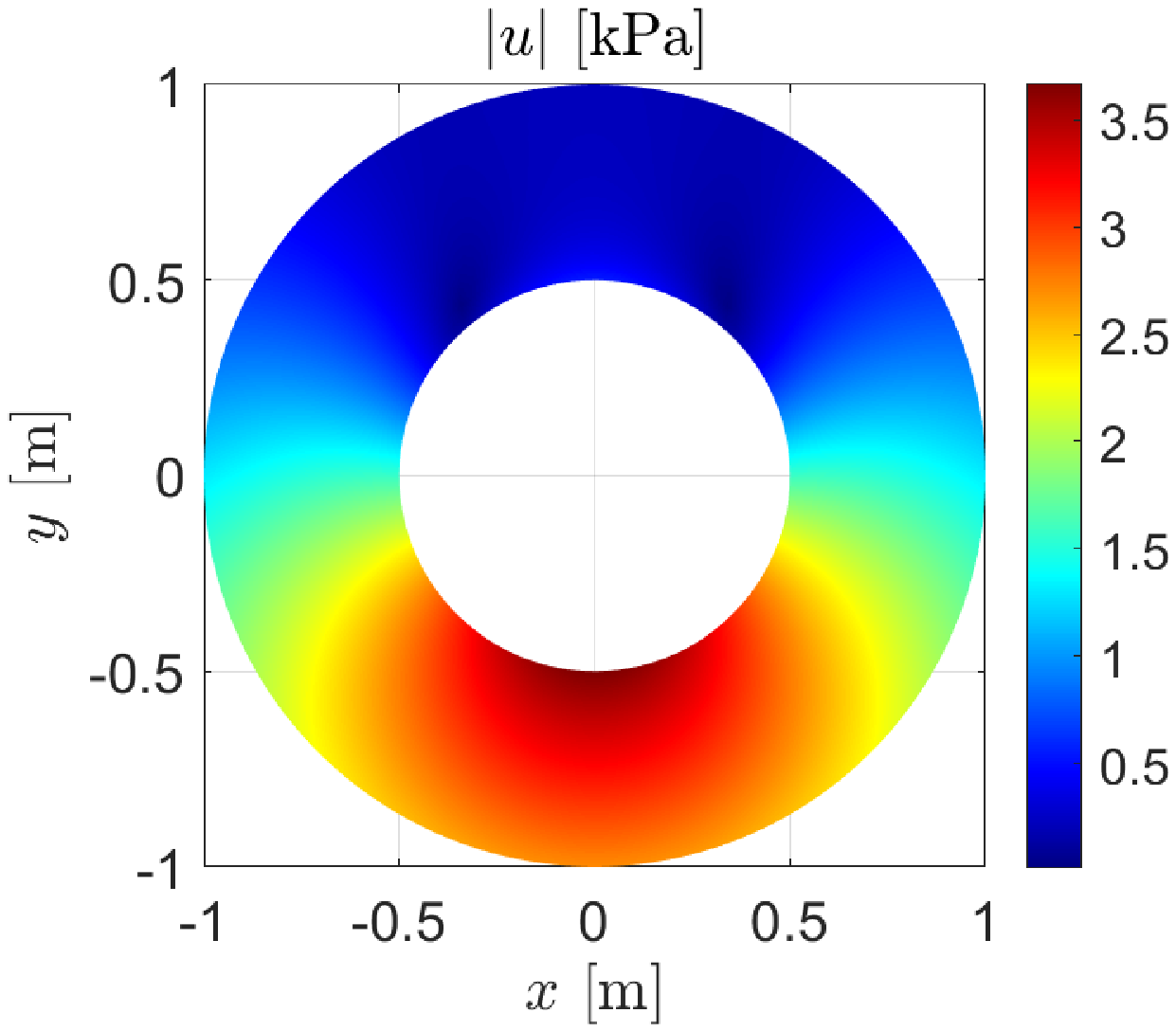}  }
    \,
    \subfloat[][]
    {   \includegraphics[trim={40 0 25 0},width=.45\textwidth]{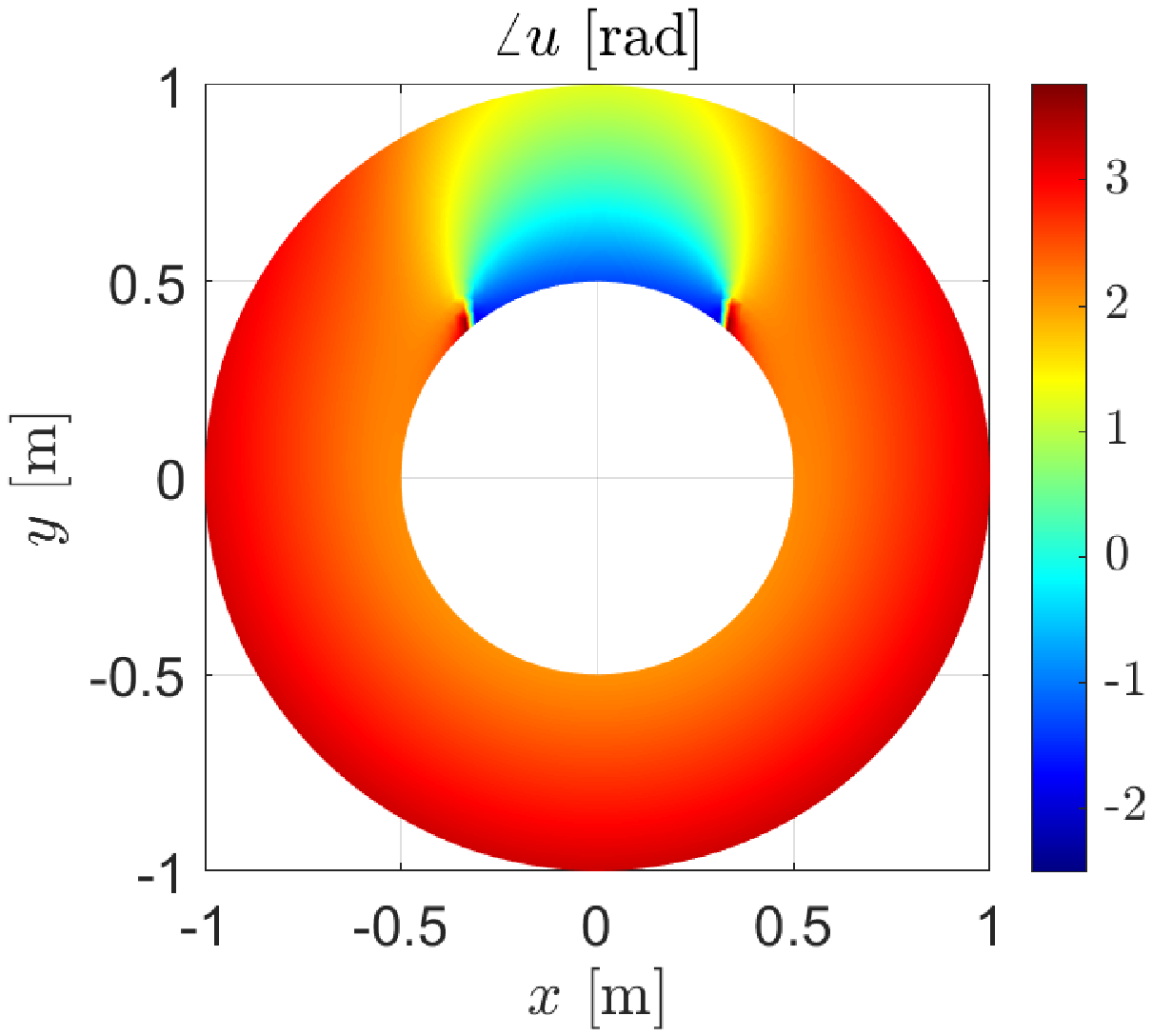}  }
    \caption{Amplitude (a) and phase (b) of the optimal pressure on $\Gamma_c$; the improvement of the oscillations $\{\bb X\}$ is summarized in table~\ref{tab:oscillation reduction}.
    }
    \label{fig:pressure optimal solution}
\end{figure}

The numerical problem is implemented in \Matlab{} thanks to the open source library redbKIT \cite{redbKIT} and the OCP is solved by inverting the linear system~\eqref{eq:pressure lin system}.
The optimal control action is shown in figure~\ref{fig:pressure optimal solution}~(a) and (b); 
table~\ref{tab:oscillation reduction} shows the motion amplitudes in terms of amplitude and phase of $\{\bb X\}$.

\begin{table}
$ \begin{array}{ccccc}
\toprule
\{\bb X\} [1\,\angle \SI{}{\radian}] & \text{Uncontrolled} & \text{Pressure-controlled} & \text{Membrane-controlled} & \text{Plate-controlled} \\
\midrule
x &  0.000 \,\angle -2.19      &    0.000 \,\angle -1.99 &  0.000\,\angle-1.33  &  0.002\,\angle-1.30  \\
y &  0.365 \,\angle -1.37  &    0.014 \,\angle -1.37 &  0.196\,\angle-0.11  &  0.572\,\angle-0.72\\
z &  1.539 \,\angle -2.99&    0.003 \,\angle -2.98  &  0.131\,\angle-3.0   &  0.102\,\angle-0.35 \\
\midrule
J                                    & \SI{1.25}{} & \SI{6.54e-3}{}  & \SI{7.05e-2}{} & \SI{2.08e-1}{} \\
\frac{1}{2}\{\bb X\}^\dag C\{\bb X\} & \SI{1.25}{} & \SI{1.07e-4}{} &   \SI{3.49e-2}{}    &   \SI{1.16e-1}{}\\
\bottomrule
\end{array} $
\caption{Comparison between body motion in case of the three different control strategies; the values describing rotations are not reported because, thanks to the symmetries of the sphere, they are almost zero and derive from numerical errors only.}
\label{tab:oscillation reduction}
\end{table}


\section{Surface tension and floating mass}
\label{sec:membrane controlled}

In case energy availability is limited, relying on a passive strategy is a reasonable alternative. However, while active controls can easily modify the wave propagation by acting as sources, a passive strategy usually relies on steering waves by locally modifying the wave speed.
The latter approach is analyzed in the following, where we design a floating device whose inertia and stiffness properties allow the symmetric control on wave speed we need.

Floating membrane dynamics has been used as a simplified model for floating ice and was proposed by \cite{keller1953water} where the effect of inertia and surface tension on propagation of water waves are considered. In particular, the loaded surface dynamics is derived by coupling the water flow with the motion of a membrane that is assumed to float on the surface, hence the vertical motion $\zeta$ of the water surface is assumed to be the same of the membrane displacement: 
\begin{align*}
    \nabla\cdot(T\nabla\zeta) &= g\rho\zeta + \rho\Phi_t +m\zeta_{tt} +mg &\text{in }\Gamma_c,
\end{align*}
where $T,m\colon\Gamma_c\to\R$ are the space varying surface tension and surface mass inertia, respectively. The left hand side accounts for the elastic restoring force, while the right hand side is the sum of the hydrostatic and hydrodynamic pressures, the inertia of the floating membrane, and its static load. We assume that the device cannot detach from the water surface;
by virtue of this, the internal forces of the membrane can be modeled as a surface tension.\\
The constant gravity force $mg$ causes a static displacement so it can be ignored while moving to the frequency domain:
\begin{align}\label{eq:membrane dynamics}
    \nabla\cdot(T\nabla\eta) &= (g\rho-\omega^2m)\eta + j\omega\rho\phi &\text{in }\Gamma_c,
\end{align}
According to this setup, the equilibrium function is $E=E(\phi,\eta,T,m)$ and takes the form of an elliptic PDE whose domain is the surface $\Gamma_c$. In the following, we control the dynamics \eqref{eq:membrane dynamics} along with \eqref{eq:potential dynamics} and \eqref{eq:body dynamics2} by acting on the control functions $T$ and $m$.

For a complete definition of the state $\eta$, we need to impose also appropriate boundary conditions on $\partial\Gamma_c$.
Since the membrane border is subject to a null shear force, it is free to oscillate and a homogeneous Neumann condition holds:
\begin{align}\label{eq:membrane BC}
    \eta_n&=0  &\text{on }\partial\Gamma_c.
\end{align}
Please note that this model assumes a horizontal equilibrium on the membrane boundary $\partial\Gamma_c$ which is hardly achievable in reality because no external forces act on it; in addition, the internal equilibrium of the surface must be guaranteed even if $T$ is space dependent..


We now define two control functions $u\in\U_{ad}$ and $v\in\V_{ad}$ that act on $T$ and $m$ such that
\begin{align*}
    T &= u,
    &
    1-\frac{\omega^2m}{g\rho} &= v.
\end{align*}
$u$ and $v$ are constrained since $T$ and $m$ are positive and $m$ must be upper limited for guaranteeing the buoyancy of the membrane. This second constraint depends on the density of the control device $\rho_c$ and its thickness $\delta$ according to the formula
$
    m = \rho_c \delta,
$
thus, in theory, no surface mass density $m$ can prevent the membrane to float if a proper thickness is chosen. However, we assume the limiting case to be ${\omega^2m = g\rho}$, that corresponds to the resonance condition of a floating object of surface mass $m$; more loaded surfaces experience a vanishing of waves \cite{keller1953water}.
Then, we define the admissible control spaces as:
\begin{equation}\label{eq:U_ad e V_ad}
\begin{aligned}
    \U_{ad} &= \{u\in\U,\, \varepsilon\leq u(\bb x) \quad\forall\bb x\in\Gamma_c\},
    \\
    \V_{ad} &= \{v\in\U,\, \varepsilon\leq v(\bb x) \leq 1-\varepsilon\quad\forall\bb x\in\Gamma_c\},
\end{aligned}
\end{equation}
where $\varepsilon$ is chosen to be \SI{1e-6}{} for avoid numerical issues in case of null control actions.\\
Note that in this case $\U=H^1(\Gamma_c,\R)\cap L^{\infty}(\Gamma_c,\R)$ is a space of real-valued functions, normed by the measure $\norm{u}^2_{\Gamma}= \int_\Gamma \alpha_u u^2+\beta_u|\nabla u|^2\,d\Gamma$, and that $\U_{ad}$, $\V_{ad}$ are convex sets.

Since there are two control actions, we redefine $J$ as
\begin{align}\label{eq:cost2}
    J\big(\{\bb X\},u,v\big) \coloneqq \frac{1}{2}\{\bb X\}^\dag C\{\bb X\} +
    \frac{1}{2}\norm{u}_\U^2 + \frac{1}{2}\norm{v}_\U^2.
\end{align}

Summing up, we aim at minimizing \eqref{eq:cost2} such that the three coupled dynamics \eqref{eq:potential dynamics}, \eqref{eq:body dynamics2} and \eqref{eq:membrane dynamics} with \eqref{eq:membrane BC} are satisfied, with $u\in\U_{ad}$, $v\in\V_{ad}$. Adopting the same strategy as before, we define the Lagrangian $\M\colon\W\times\W'\times\C^6\times\C^6\times\Y\times\Y'\times\U\times\U\to\R$ as
\begin{multline}\label{eq:membrane lagrangian}
    \M\big(\phi^s,\lambda,\{\bb X\},\{\bb Y\},\eta,\mu, u,v\big) \coloneqq
    J
    +\Re\Big\{
    \int_\Omega{\nabla\bar\lambda\cdot\nabla\phi^s}
    -\frac{\omega^2}{g}\int_{\Gamma_f}{\bar\lambda\phi^s}
    +\int_{\Gamma_c}{\bar\lambda(\phi^i_n-j\omega\eta)}
    \\
    -j\omega\int_{\Gamma_g}{\bar\lambda\{\bb n\}^\dag}\,\{\bb X\}
    +\int_{\Gamma_g}{\bar\lambda\phi^i_n}
    + \alpha\int_{\Gamma_e}{\bar\lambda\phi^s}
    \\
    +\{\bb Y\}^\dag
    \Big[\big(K-\omega^2M\big)\{\bb X\} +j\omega\rho\int_{\Gamma_g}(\phi^s+\phi^i)\{\bb n\}\Big]
    \\
    + \int_{\Gamma_c}{u\nabla\bar\mu\cdot\nabla\eta + g\rho v\bar\mu\eta + j\omega\rho\bar\mu(\phi^s+\phi^i)}
    \Big\},
\end{multline}
where $\eta\in\Y$, $\Y=H^1(\Gamma_c)$, and $\mu\in\Y'$ the adjoint state of $\eta$, $\Y'$ its adjoint space. The first order necessary conditions are obtained  by taking the derivatives of $\M$ with respect to the states $\phi^s$, $\eta$ and $\{\bb X\}$ and controls $u$, $v$. For the sake of brevity, the strong formulations are reported only:
\begin{align*}
\text{adjoint:} \,\, &
    \begin{aligned}
    &\begin{dcases}
    -\Delta\lambda = 0 &\text{in }\omega
    \\
    \lambda_n = 0 &\text{on }\Gamma_r
    \\
    \lambda_n - \frac{\omega^2}{g}\lambda= 0 &\text{on }\Gamma_f
    \\
    \lambda_n = j\omega\rho\mu &\text{on } \Gamma_c
    \\
    \lambda_n = j\omega\rho \{\bb n\}^\top\{\bb Y\} &\text{on }\Gamma_g
    \\
    \lambda_n +\bar\alpha\lambda = 0 &\text{on }\Gamma_e
    \end{dcases}
    \\
    &(K-\omega^2M)^\top\{\bb Y\} = -j\omega\int_{\Gamma_g}{\lambda\{\bb n\}} - C^\top\{\bb X\}
    \\
    &\begin{dcases}
    -\nabla\cdot(u\nabla\mu) + g\rho v\mu +j\omega\lambda = 0& \text{on }\Gamma_c
    \\
    \mu_n = 0 &\text{on }\partial\Gamma_c
    \end{dcases}
\end{aligned}
\\[5pt] \text{control:}    \,\, &
\begin{aligned}
    \LL'_u[\psi-u^*] &= \int_{\Gamma_c}\beta_u \nabla u^*\cdot\nabla(\psi-u^*) + \big[\alpha_u u^* + \Re( \nabla\bar\mu^*\cdot\nabla\eta^* )\big](\psi-u^*)  \geq 0 & \forall\psi\in\U_{ad}
    \\
    \LL'_v[\psi-v^*] &= \int_{\Gamma_c} \beta_v\nabla v^*\cdot\nabla(\psi-v^*) +\big[\alpha_v v^* + g\rho\Re( \bar\mu^*\eta^* )\big](\psi-v^*)\geq 0 & \forall\psi\in\V_{ad}
\end{aligned}
\end{align*}
where $(\phi^*,\lambda^*,\eta^*,\mu^*,u^*,v^*)$ is the optimal solution of the control problem.
Note that this is a nonlinear OCP, thus the solution must be obtained by means of an iterative approach.

Relying on the same discretization as before, we obtain the following optimality conditions
\begin{equation}\label{eq:membrane discrete}
\begin{aligned}
    \text{States:} & \qquad
    \begin{aligned}
        &\left[ A-\frac{\omega^2}{g}C_f+\alpha C_e \right] \bs\phi = j\omega K_g^\top\{\bb X\} + j\omega D_c \bs\eta + \bb f_g + \bb f_c
        \\
        &(K-\omega^2M)\{\bb X\}=-j\omega\rho (K_g\bs\phi + \bb g)
        \\
        & (\mathbb A\bb u + \rho g \mathbb B\bb v)\,\bs\eta = -j\omega\rho D_c^\top \bs\phi
    \end{aligned}
    \\[5pt]
    \text{Adjoint states:} & \qquad
    \begin{aligned}
        &\left[ A - \frac{\omega^2}{g}C_f +\bar\alpha C_e \right]\bs\lambda = j\omega\rho K_g^\top\{\bb Y\} + j\omega\rho D_c\bs\mu
        \\
        &(K-\omega^2M)^\top\{\bb Y\}=-j\omega K_g\bs\lambda - C^\top\{\bb X\}
        \\
        & (\mathbb A\bb u + \rho g \mathbb B\bb v)\,\bs\mu = -j\omega D_c^\top \bs\phi
    \end{aligned}
    \\[5pt]
    \text{Control gradients:} &\qquad
    \begin{aligned}
        &\nabla_{\bb u} J = \left(\alpha_u E_c +\beta_u A_c\right)\bb u + \Re( \bar{\bs\mu}\mathbb A\,\bs\eta)
        \\
        &\nabla_{\bb v} J = \left(\alpha_v E_c+\beta_v A_c\right)\bb v + g\rho\,\Re(\bar{\bs\mu}\mathbb B\,\bs\eta)
    \end{aligned}
\end{aligned}
\end{equation}
where the third order tensors $\mathbb A,\mathbb B\in\R^{n\times n\times l}$ are defined as
\begin{align*}
    (\mathbb A)_{ijk} \coloneqq \int_{\Gamma_c} \nabla\phi_i\cdot\nabla\phi_j\,\psi_k\,d\Gamma,
    &&
    (\mathbb B)_{ijk} \coloneqq \int_{\Gamma_c} \phi_i\phi_j\psi_k\,d\Gamma,
\end{align*}
and the products $\mathbb A\bb u$ and $\bs\mu\mathbb A$ are matrices defined as
\begin{align*}
    (\mathbb A\bb u)_{ij} \coloneqq \sum_{k=1}^l (\mathbb A)_{ijk}\,u_k,
    &&
    (\bs\mu\mathbb A\bb)_{kj} \coloneqq \sum_{i=1}^n \mu_i\,(\mathbb A)_{ijk}.
\end{align*}

\begin{figure}
    \centering
    \subfloat[][]
    {   \includegraphics[trim={14 -2 10 0},width=.325\textwidth]{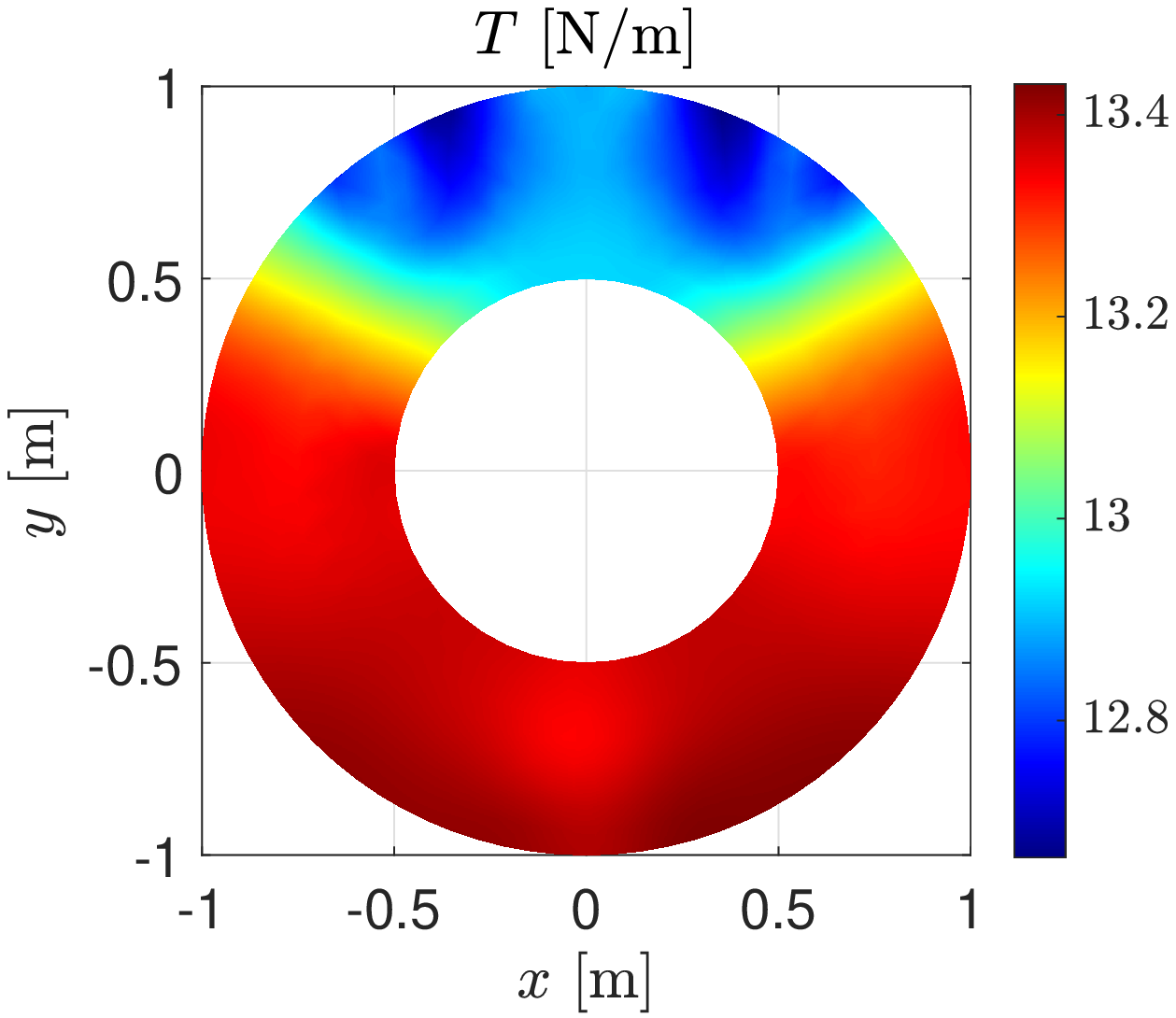}  }
    \,
    \subfloat[][]
    {   \includegraphics[trim={10 0 10 0},width=.325\textwidth]{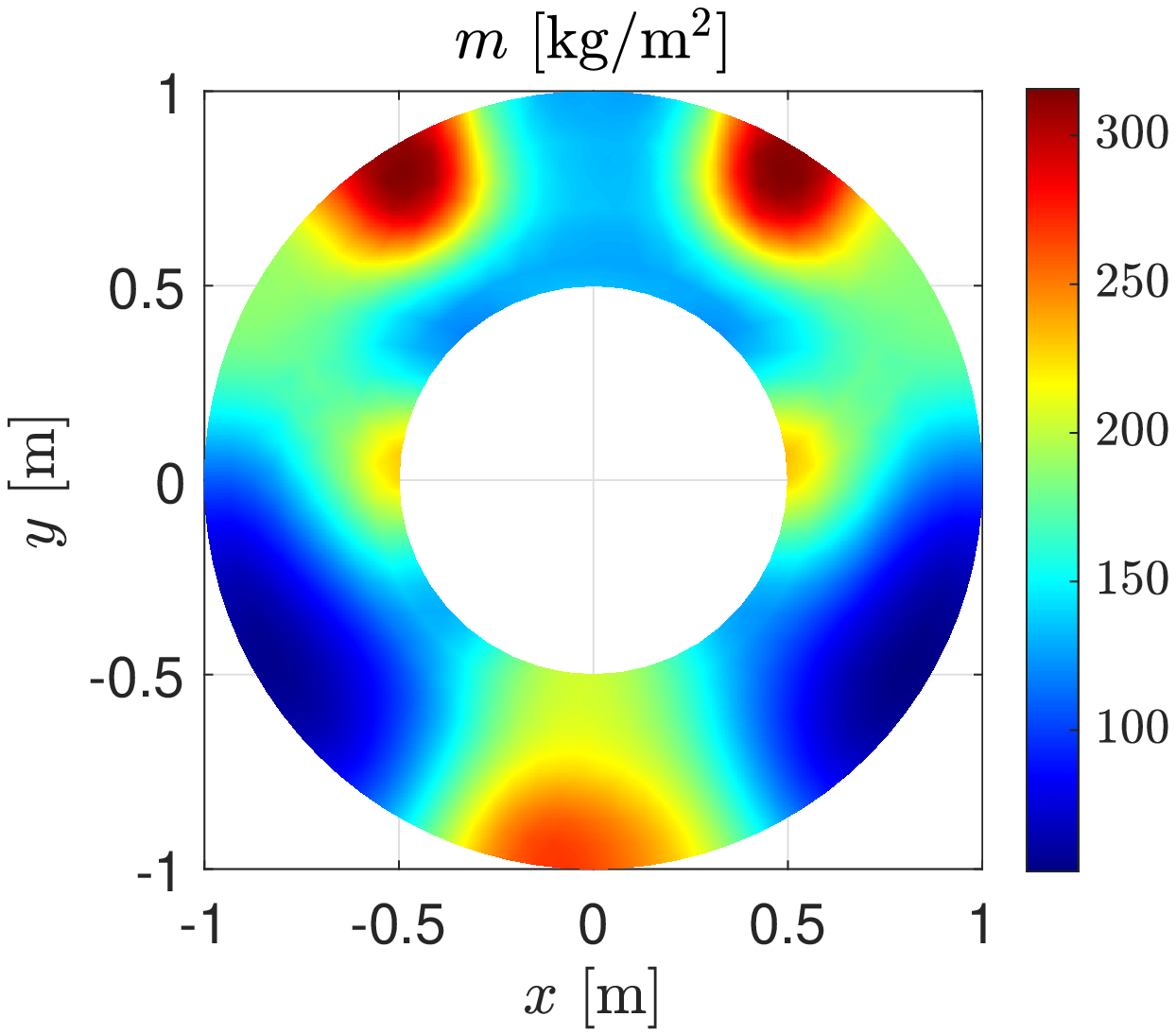}  }
    \,
    \subfloat[][]
    {   \includegraphics[trim={5 -2 15 0},width=.28\textwidth]{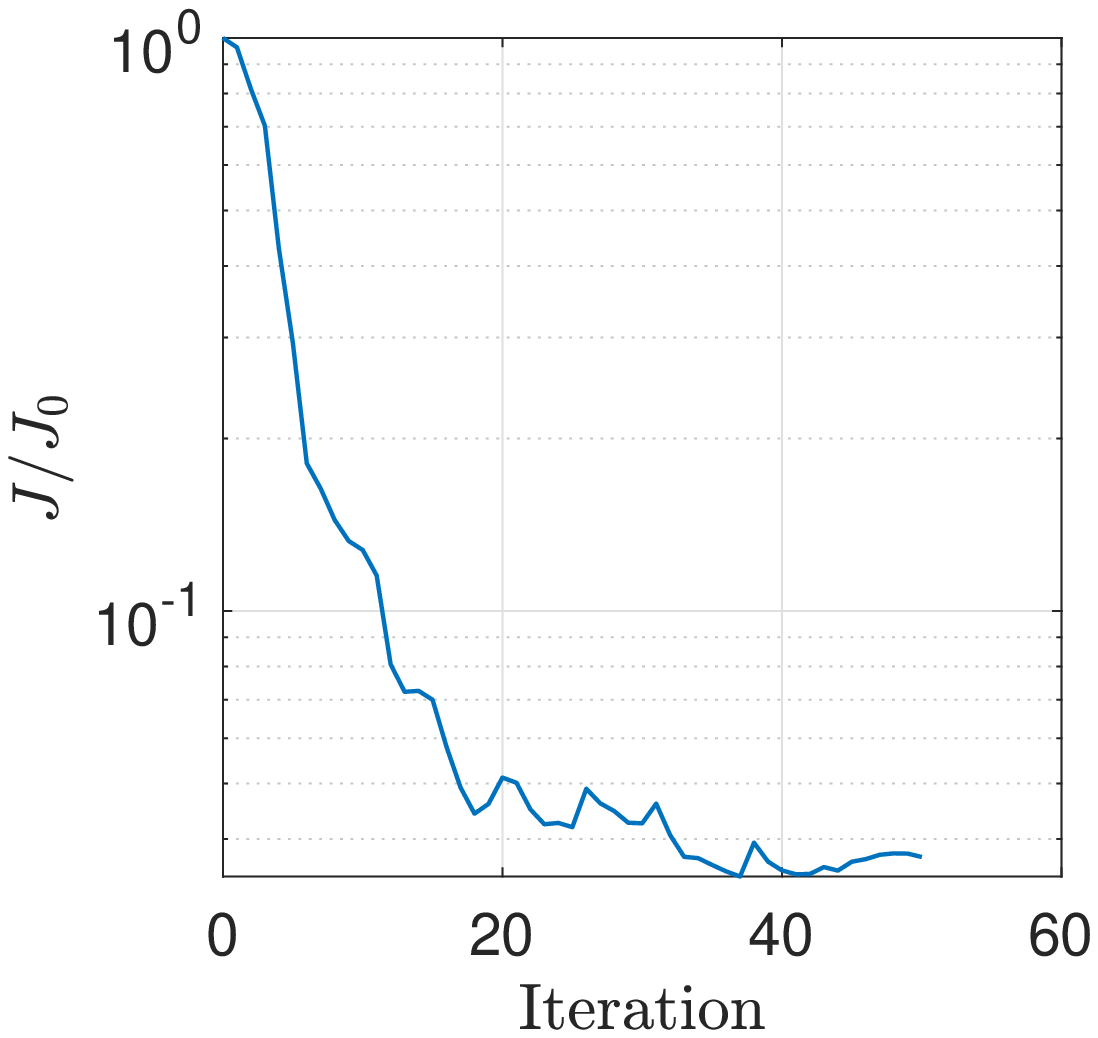}}
    \caption{Membrane properties in terms of optimal controls u (a) and v (b); (c) objective over iteration.  The improvement of the oscillations $\{\bb X\}$ is summarized in table~\ref{tab:oscillation reduction}.
    }
    \label{fig:membrane optimal solution}
\end{figure}

We consider the same setup adopted for the pressure-controlled case, but a model of the floating membrane is added on the surface $\Gamma_c$. We set $\alpha_u=\alpha_v=\SI{1e-4}{}$ and $\beta_u=\beta_v=\SI{4e-2}{}$; then, we solve the nonlinear OCP relying on the \Matlab{} interface for IPOPT by~\cite{wachter2006implementation,mexIPOPT}, a software package for large-scale nonlinear optimization. The resulting optimal $T$ and $m$ are shown in figure~\ref{fig:membrane optimal solution} along with the optimization convergence in terms of cost minimization.

\section{Floating plate}
\label{sec:plate controlled}

The control strategy just proposed seems effective on reducing the floating body motion; however a device with variable surface tension can be of difficult realization due to the problems of internal equilibrium discussed above, then a more sophisticated control can be useful in practice. In particular, the elastic restoring force given by surface tension must be replaced by another one.  \\
Through this section, the control surface $\Gamma_c$ is loaded by a floating elastic thin plate. The model we adopt is similar to the one proposed in literature for studying the effect of floating ice on water waves, see e.g.\ \cite{balmforth1999ocean,fox1990reflection}. In this setting, the plate stiffness acts as a restoring force thus increasing the wave speed, conversely the load inertia slows it down. So, similarly to the previous case, we expect to control the water flow by deriving stiffness and inertia properties of the plate.

Let us consider the plate equilibrium in weak form \cite{timoshenko1959theory}:
\begin{align}\label{eq:plate weak}
    \int_{\Gamma_c} {B\big[ (1-\nu)D^2 \eta : D^2  v + \nu\Delta\eta\Delta  v\big]\,d\Gamma} &= \int_{\Gamma_c} {f v \,d\Gamma}& \forall v\in \Y,
\end{align}
where $B\colon\Gamma_c\to\R$ is the local flexural stiffness, $\nu$ is the Poisson ratio of the plate material,  $f$ is a distributed load pointing upwards. $D^2\eta$ is the hessian matrix of $\eta$ and $A:B$ indicates the Frobenius scalar product between matrices $A$ and $B$, i.e.\
$
    A:B \coloneqq \sum_{i,j=1}^2 (A)_{ij}(B)_{ij}.
$
Note that in this case we need to require that $\eta\in\Y=H^2(\Gamma_c)$ for a correct definition of the weak formulation \eqref{eq:plate weak}.
By following a procedure similar to the one in \cite{timoshenko1959theory}, we can derive the strong form in case of space dependent $B$:
\begin{equation}\label{eq:plate dynamics}
\begin{dcases}
\Delta(B\Delta \eta) + (1-\nu)\left(2B_{xy}\eta_{xy} - B_{xx}\eta_{yy} - B_{yy}\eta_{xx}\right) - f  = 0  & \text{in } \Gamma_c
\\
\Delta \eta -(1-\nu)\eta_{\tau\tau} = 0 & \text{on } \partial\Gamma_c
\\
(B\Delta \eta)_n + (1-\nu)(2 B_\tau \eta_{n\tau} - B_n \eta_{n\tau} + B\eta_{n\tau\tau}) = 0 & \text{on } \partial\Gamma_c
\end{dcases},
\end{equation}
 where the subscript $\tau$ stands for the tangential derivative on the boundary.
Note that the two conditions on $\partial\Gamma_c$ holds because we impose the plate edges to be free, i.e.\ null torque and null shear.
The force $f$ acting on the plate is given by the inertia of the plate itself and the water pressure obtained from \eqref{eq:lin Bernoulli}, then
\begin{equation*}
    f = -m\zeta_{tt} - mg - \rho g\zeta - \rho\Phi_t,
\end{equation*}
and, in frequency domain,
\begin{equation*}
    f = (m\omega^2 - \rho g)\eta - j\omega\rho\phi.
\end{equation*}

The control problem reads: find $u$, $v$ that minimize \eqref{eq:cost2} where the three coupled dynamics \eqref{eq:potential dynamics}, \eqref{eq:body dynamics2} and \eqref{eq:plate dynamics} are satisfied, with $u\in\U_{ad}$, $v\in\V_{ad}$ where
\begin{align*}
    B &= u
    &
    1-\frac{\omega^2m}{g\rho} &= v;
\end{align*}
and $\U_{ad}$, $\V_{ad}$ are defined by \eqref{eq:U_ad e V_ad}.
Again, we apply the Lagrange's multipliers method for obtaining the first order necessary conditions, similarly to \eqref{eq:membrane lagrangian} we define
$\N\colon\W\times\W'\times\C^6\times\C^6\times\Y\times\Y'\times\U\times\U\to\R$ as
\begin{multline}\label{eq:plate lagrangian}
    \N\big(\phi^s,\lambda,\{\bb X\},\{\bb Y\},\eta,\mu, u,v\big) \coloneqq
    J
    +\Re\Big\{
    \int_\Omega{\nabla\bar\lambda\cdot\nabla\phi^s}
    -\frac{\omega^2}{g}\int_{\Gamma_f}{\bar\lambda\phi^s}
    +\int_{\Gamma_c}{\bar\lambda(\phi^i_n-j\omega\eta)}
    \\
    -j\omega\int_{\Gamma_g}{\bar\lambda\{\bb n\}^\dag}\,\{\bb X\}
    +\int_{\Gamma_g}{\bar\lambda\phi^i_n}
    + \alpha\int_{\Gamma_e}{\bar\lambda\phi^s}
    +\{\bb Y\}^\dag
    \Big[\big(K-\omega^2M\big)\{\bb X\} +j\omega\rho\int_{\Gamma_g}(\phi^s+\phi^i)\{\bb n\}\Big]
    \\
    +\int_{\Gamma_c} {u\big[ (1-\nu)D^2 \eta : D^2\bar\mu + \nu\Delta\eta\Delta\bar\mu \big]} + \int_{\Gamma_c} {\big[g\rho v\eta + j\omega\rho(\phi^s+\phi^i)\big] \bar\mu}
    \Big\}.
\end{multline}
Following the same procedure as before, we obtain the adjoints equations in the strong form
\begin{equation}\label{eq:plate adjoint}
\begin{gathered}
    \begin{dcases}
    -\Delta\lambda = 0 &\text{in }\omega
    \\
    \lambda_n = 0 &\text{on }\Gamma_r
    \\
    \lambda_n - \frac{\omega^2}{g}\lambda= 0 &\text{on }\Gamma_f
    \\
    \lambda_n = j\omega\rho\mu &\text{on } \Gamma_c
    \\
    \lambda_n = j\omega\rho \{\bb n\}^\top\{\bb Y\} &\text{on }\Gamma_g
    \\
    \lambda_n +\bar\alpha\lambda = 0 &\text{on }\Gamma_e
    \end{dcases}
    \\[5pt]
    (K-\omega^2M)^\top\{\bb Y\} = -j\omega\int_{\Gamma_g}{\lambda\{\bb n\}\,d\Gamma} - C^\top\{\bb X\}
    \\[4pt]
    \begin{dcases}
    \Delta(u\Delta \mu) + (1-\nu)\left(2u_{xy}\mu_{xy} - u_{xx}\mu_{yy} - u_{yy}\mu_{xx}\right) +\rho g v\mu +j\omega\lambda  = 0  & \text{in } \Gamma_c
    \\
    \Delta \mu -(1-\nu)\mu_{\tau\tau} = 0 & \text{on } \partial\Gamma_c
    \\
    (u\Delta \mu)_n + (1-\nu)(2 u_\tau \mu_{n\tau} - u_n \mu_{n\tau} + u\mu_{n\tau\tau}) = 0 & \text{on } \partial\Gamma_c
    \end{dcases}
\end{gathered}\end{equation}
and the control gradients
\begin{equation*}
\begin{aligned}
    \LL'_u[\psi-u^*] &= \int_{\Gamma_c}\beta_u \nabla u^*\cdot\nabla(\psi-u^*) + \big[\alpha_u u^* + \Re( \nabla\bar\mu^*\cdot\nabla\eta^* )\big](\psi-u^*)  \geq 0 & \forall\psi\in\U_{ad}
    \\
    \LL'_v[\psi-v^*] &= \int_{\Gamma_c} \beta_v\nabla v^*\cdot\nabla(\psi-v^*) +\big[\alpha_v v^* + g\rho\,\Re( \bar\mu^*\eta^* )\big](\psi-v^*)\geq 0 & \forall\psi\in\V_{ad}
\end{aligned}\,.
\end{equation*}

Finally, the first order optimality conditions are discretized in order to obtain a numerical solution. With respect to the discrete nonlinear system \eqref{eq:membrane discrete}, the only difference we have is that the stiffness tensor $\mathbb A$, describing the elastic restoring force due to internal tension, is substituted by $\mathbb L$ that accounts for the flexural stiffness. For this reason, we do not rewrite the entire set of equations. However, it is important to highlight how $\mathbb L$ is defined, that is
\begin{equation}
(\mathbb L)_{ijk} = \int_{\Gamma_c} {\psi_k \big[ (1-\nu)D^2 \phi_i : D^2\phi_j + \nu\Delta\phi_i\Delta\phi_j\big]\,d\Gamma},
\end{equation}
that is numerical consistent only if $\W_h\subseteq H^2(\Omega)$, which is not the case for a FEM based on $\mathbb P_1$ shape functions. In order to overcome this issue while avoiding the computational effort required by standard nonconforming or mixed methods \cite{joly2005scientific}, we adopt the recovery gradient technique described in \cite{guo2018c} that allows to solve biharmonic problems relying on $C^0$ finite elements.

\begin{figure}
    \centering
    \subfloat[][]
    {   \includegraphics[trim={14 0 10 0},width=.325\textwidth]{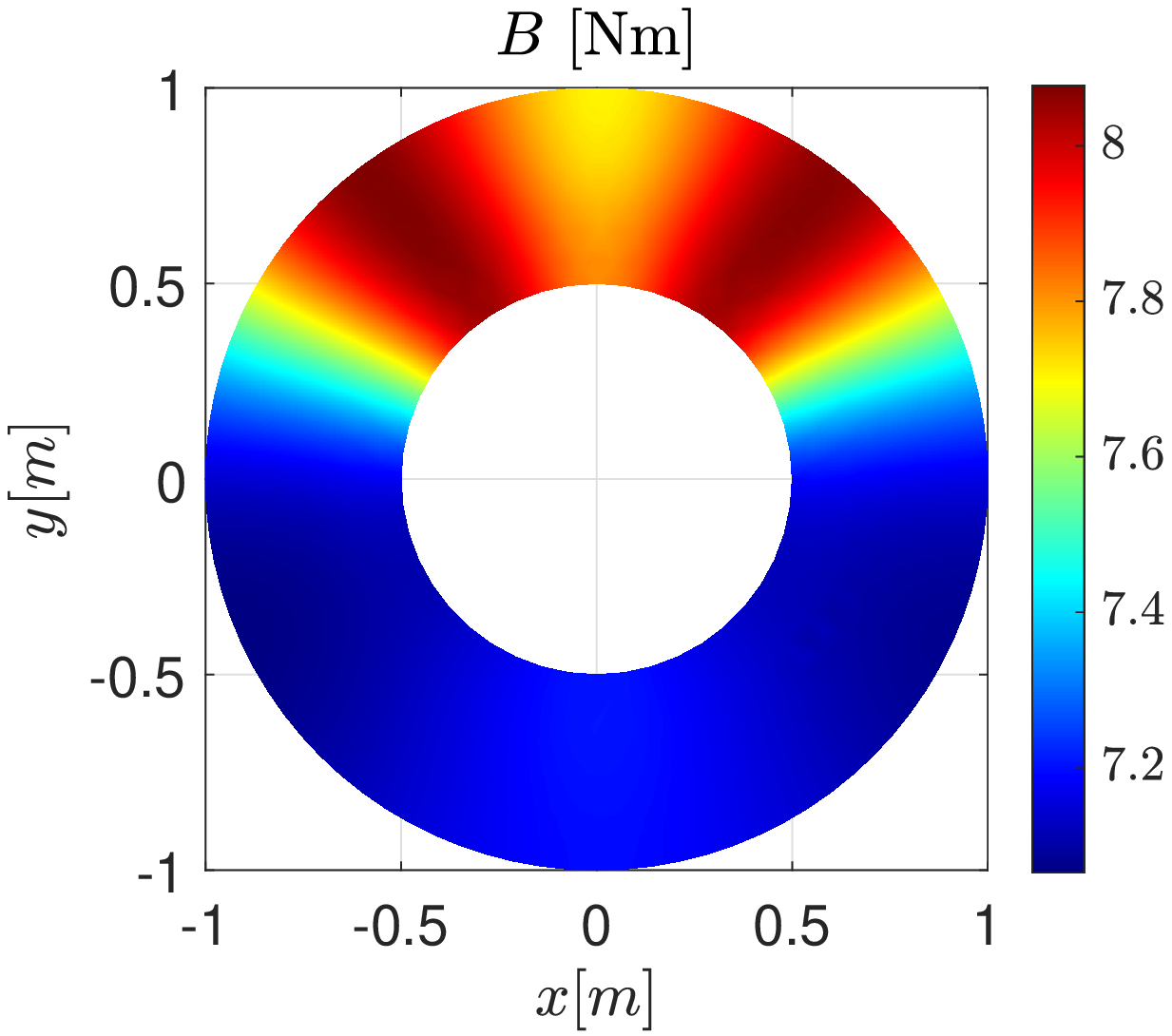}  }
    \,
    \subfloat[][]
    {   \includegraphics[trim={10 0 10 0},width=.325\textwidth]{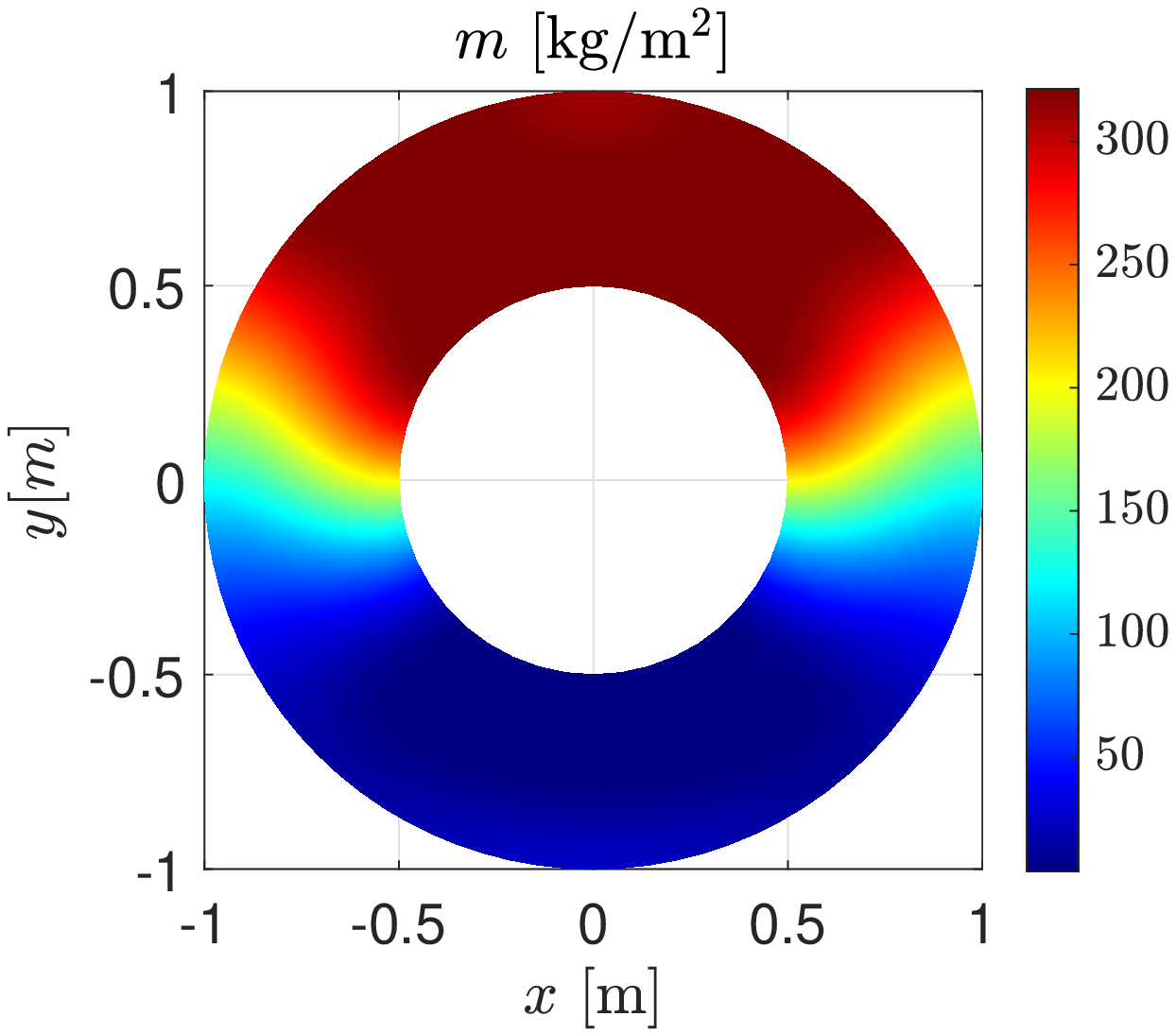}  }
    \,
    \subfloat[][]
    {   \includegraphics[trim={5 -2 15 0},width=.28\textwidth]{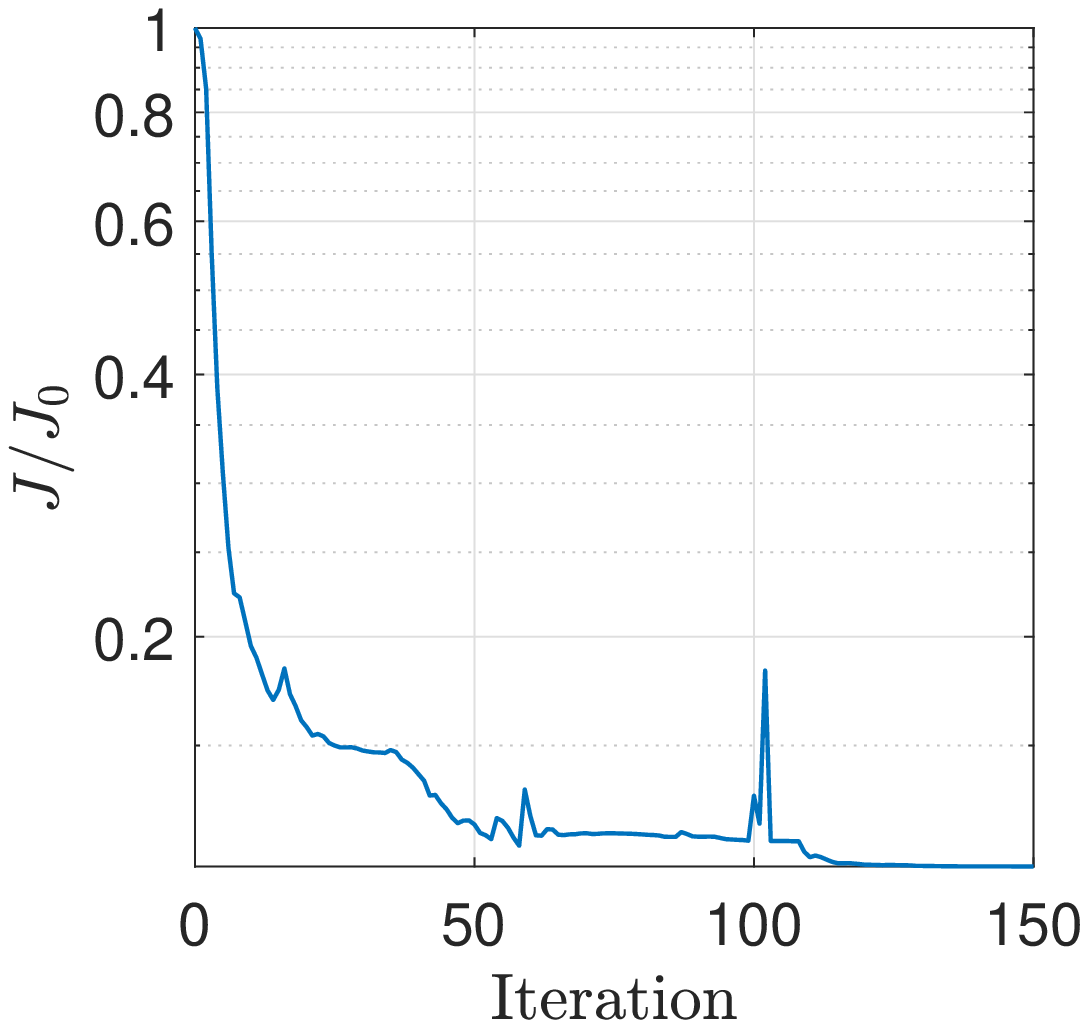}}
    \caption{Plate properties in terms of optimal controls u (a) and v (b); (c) objective over iteration. The improvement of the oscillations $\{\bb X\}$ is summarized in table~\ref{tab:oscillation reduction}.
    }
    \label{fig:plate optimal solution}
\end{figure}

The setup studied before is now controlled by means of a plate whose properties are determined by solving the corresponding nonlinear OCP. Again, the minimization is carried out making use of IPOPT software; the results in terms of optimal stiffness and inertia properties are shown in figure \ref{fig:plate optimal solution} and the vibration reductions are highlighted in table~\ref{tab:oscillation reduction}.

\section{Floating wind turbine}
\label{sec:turbine}
\begin{figure}
\centering\label{fig:turbine}
    \begin{tabular}{ccc}
        \multirow{2}{*}[55pt]{
            \subfloat[]{\label{fig:turb a}
            \includegraphics[width=0.25\hsize,trim={30 -35 20 20}]{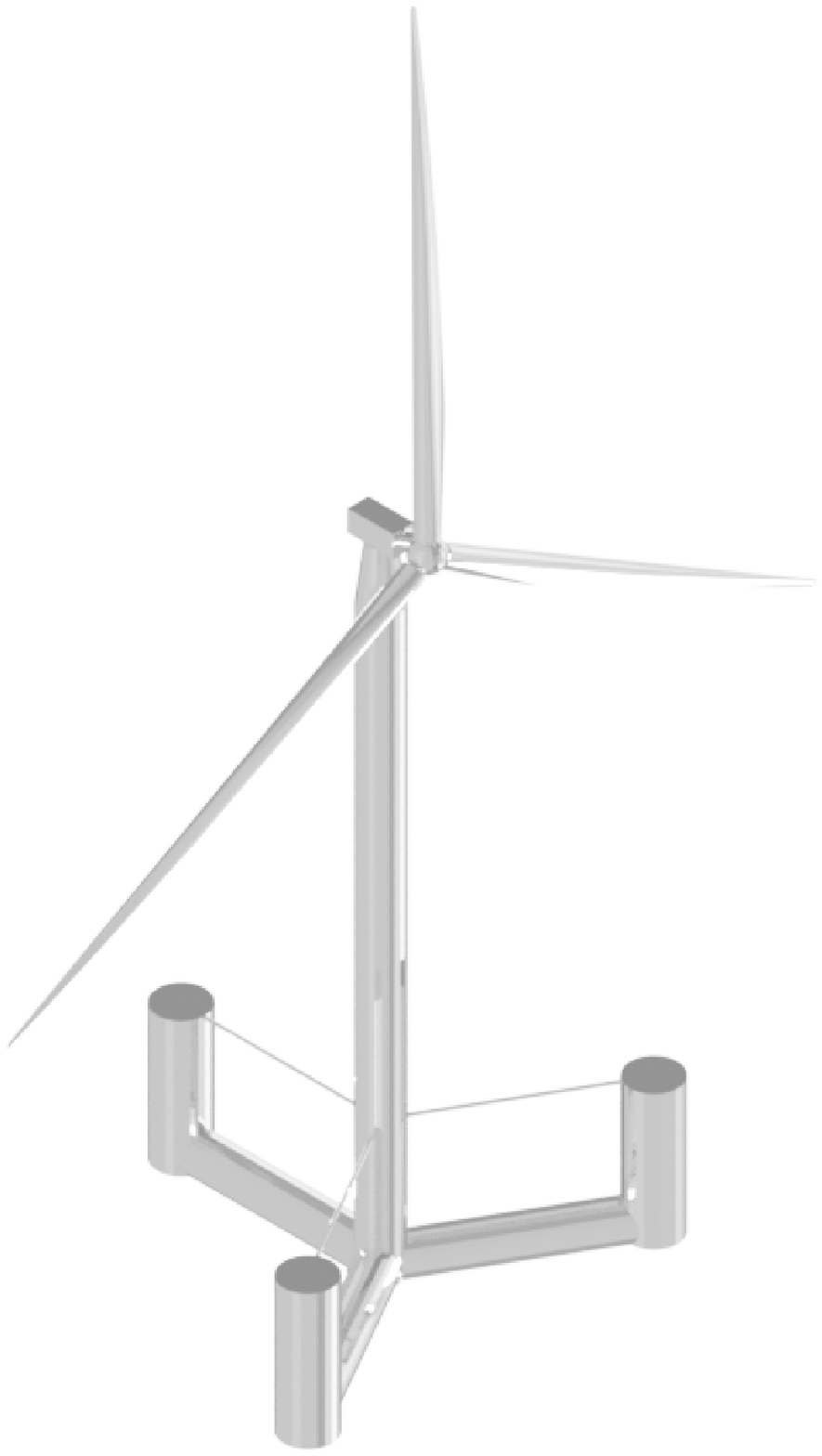}
            }
        }
        &
        \subfloat[]{\label{fig:turb b}
        \includegraphics[width=0.32\textwidth,trim={40 0 0 -10}]{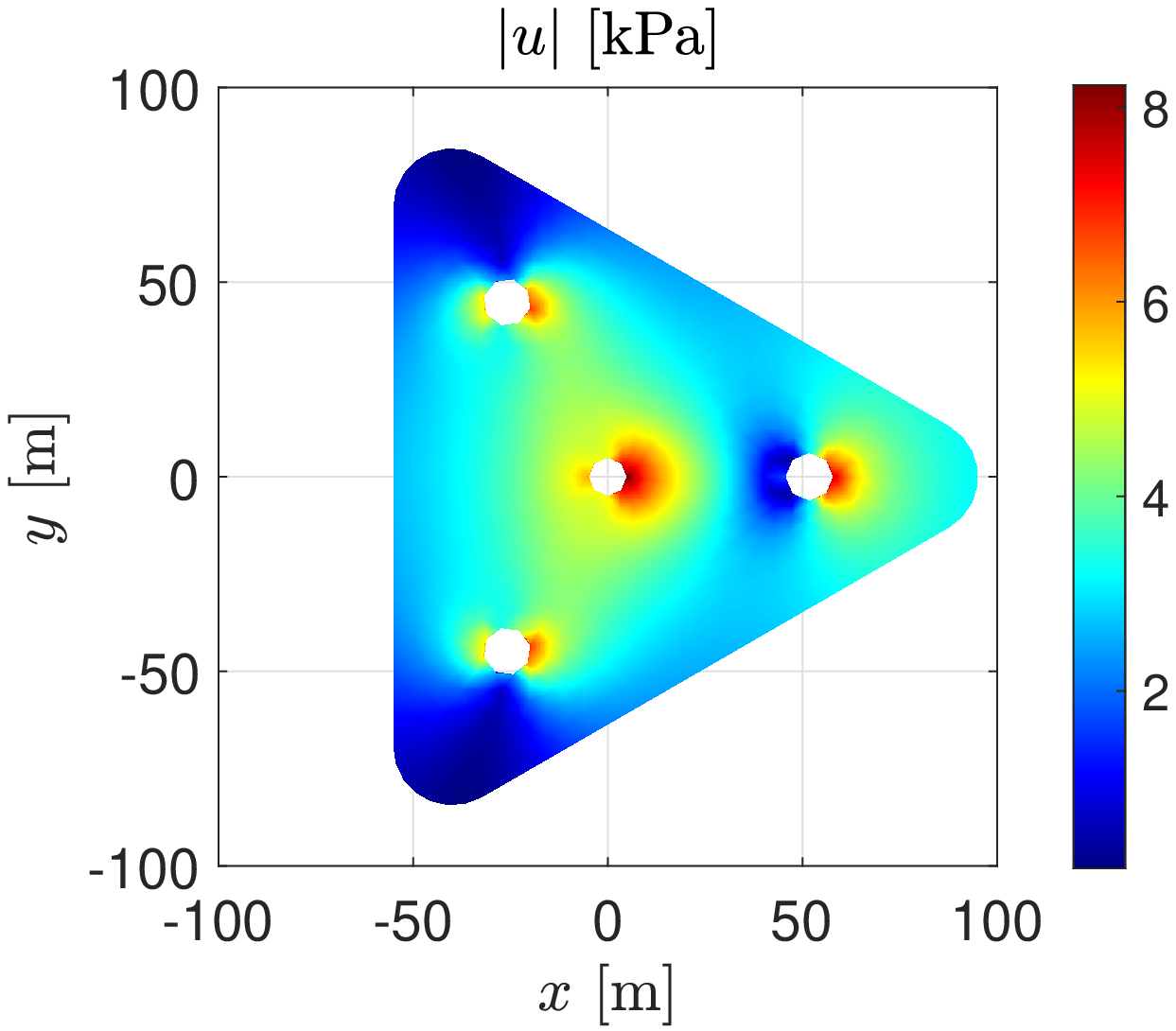}}
        &
        \subfloat[]{\label{fig:turb c}
        \includegraphics[width=0.32\textwidth,trim={20 0 20 -10}]{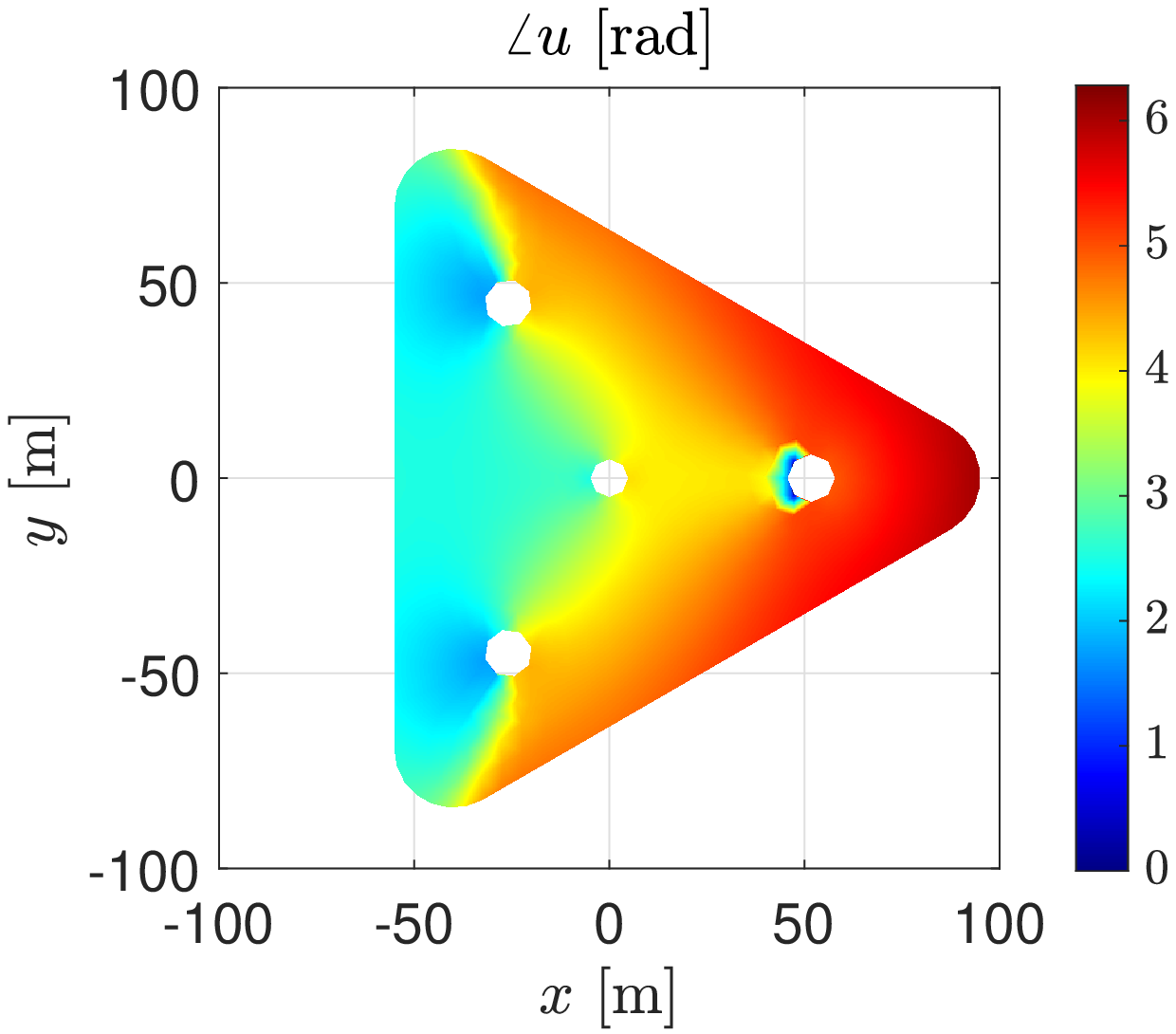}}
        \\
        &
        \subfloat[]{\label{fig:turb d}
        \includegraphics[width=0.32\textwidth,trim={40 0 14 -10}]{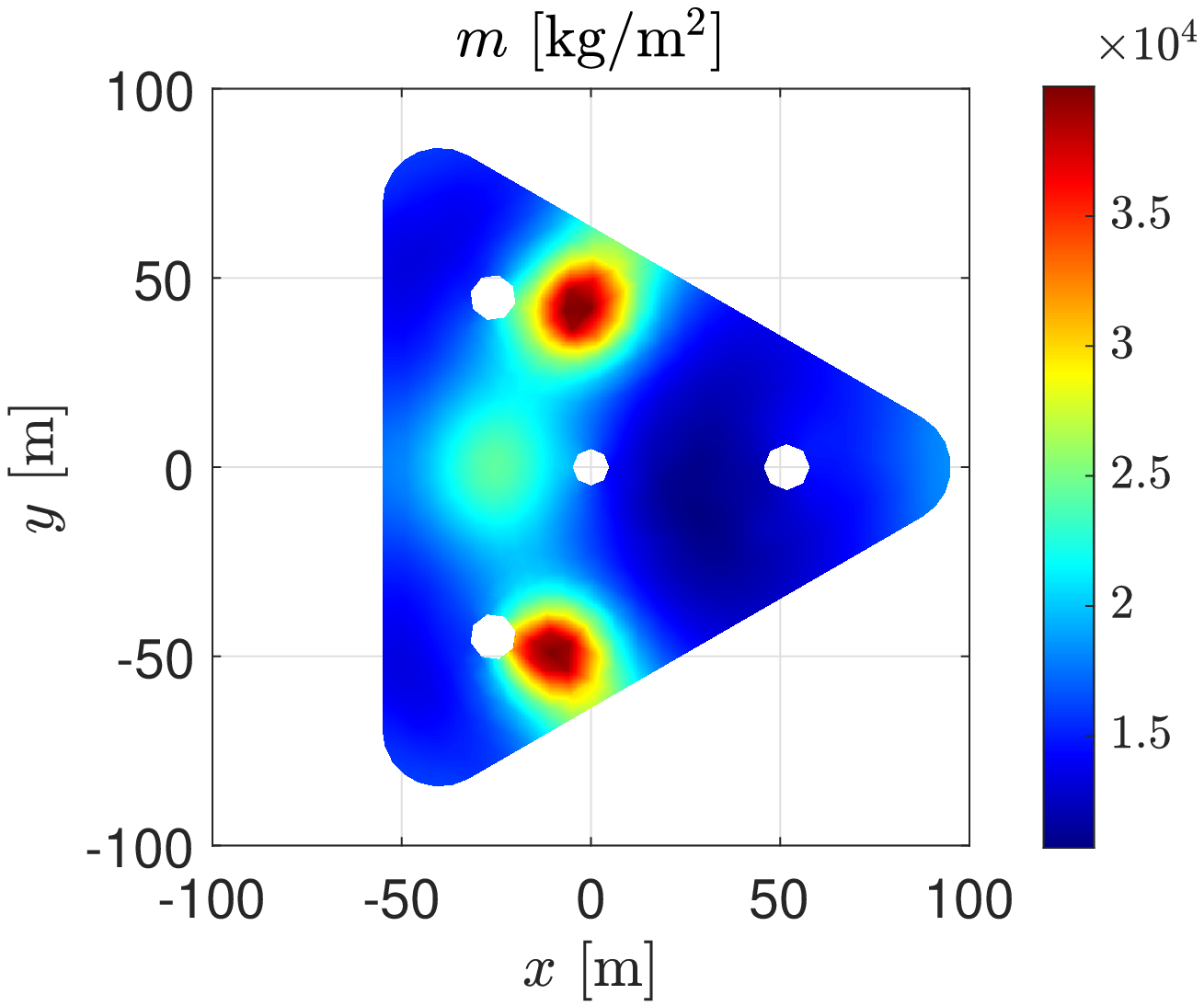}}
        &
        \subfloat[]{\label{fig:turb e}
        \includegraphics[width=0.32\textwidth,trim={20 0 23 -10}]{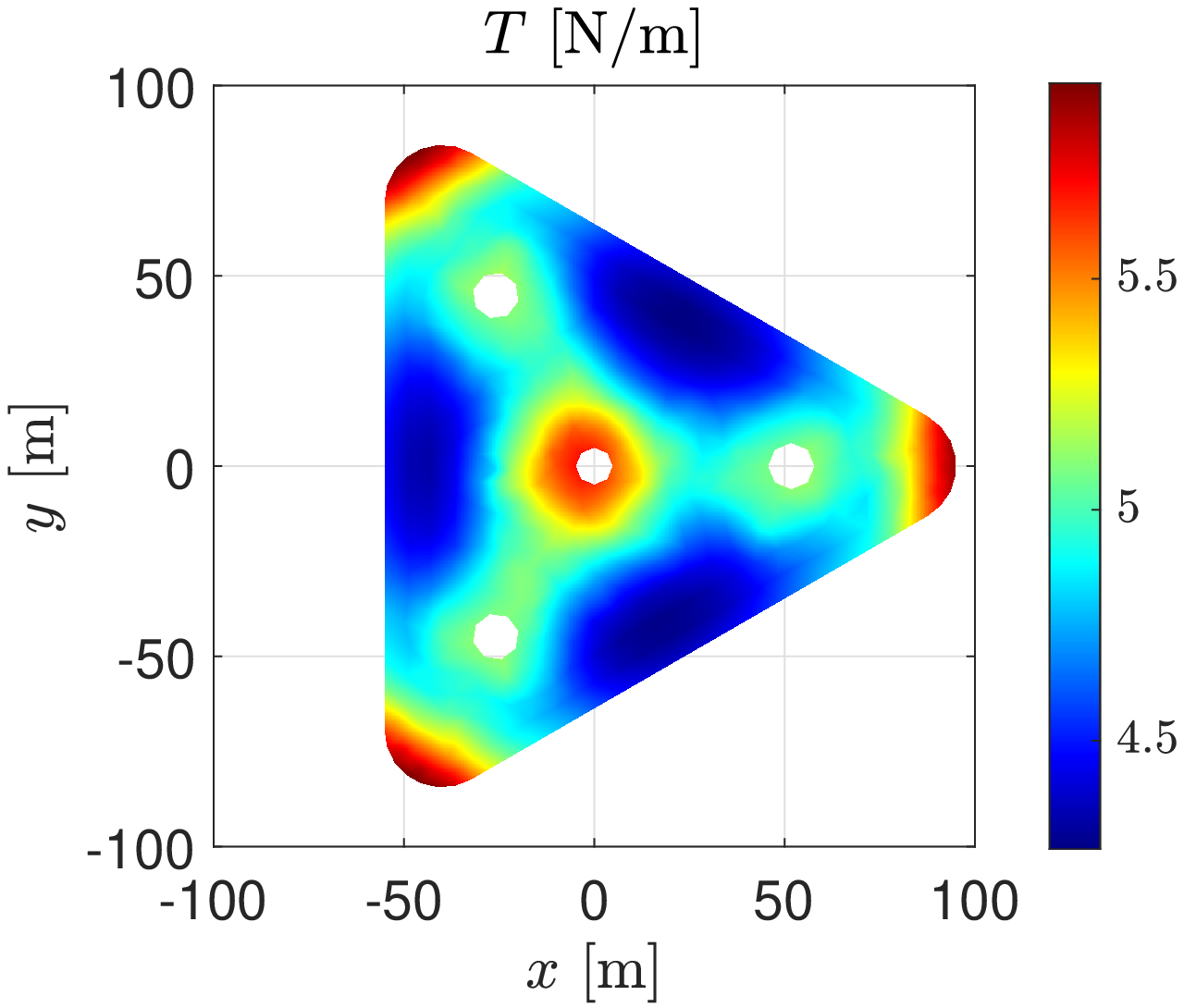}}
        \\
        \subfloat[]{\label{fig:turb f}
        \includegraphics[width=0.28\textwidth,trim={36 -5 -4  40 -10}]{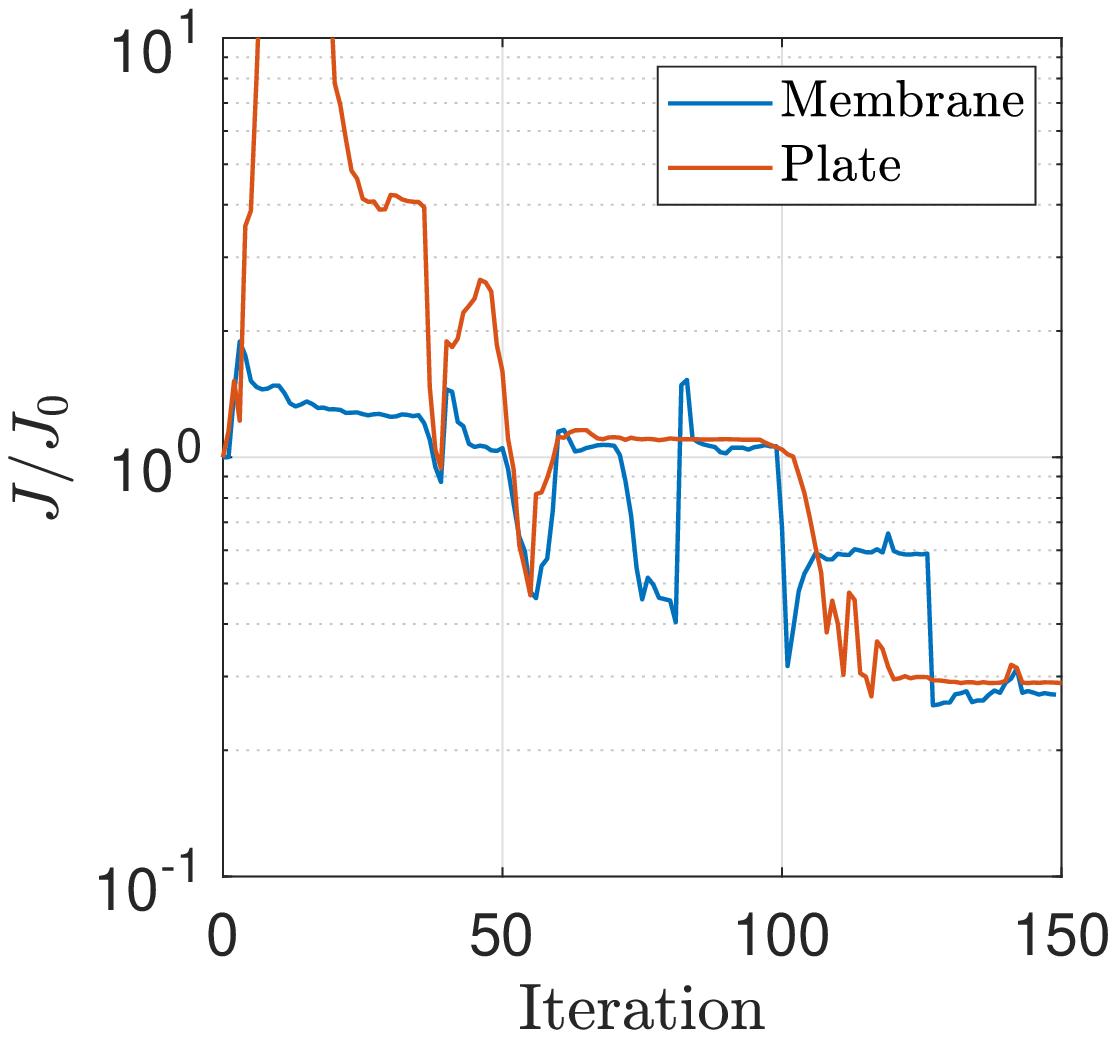}}
        &
        \subfloat[]{\label{fig:turb g}
        \includegraphics[width=0.32\textwidth,trim={40 0 10 -10}]{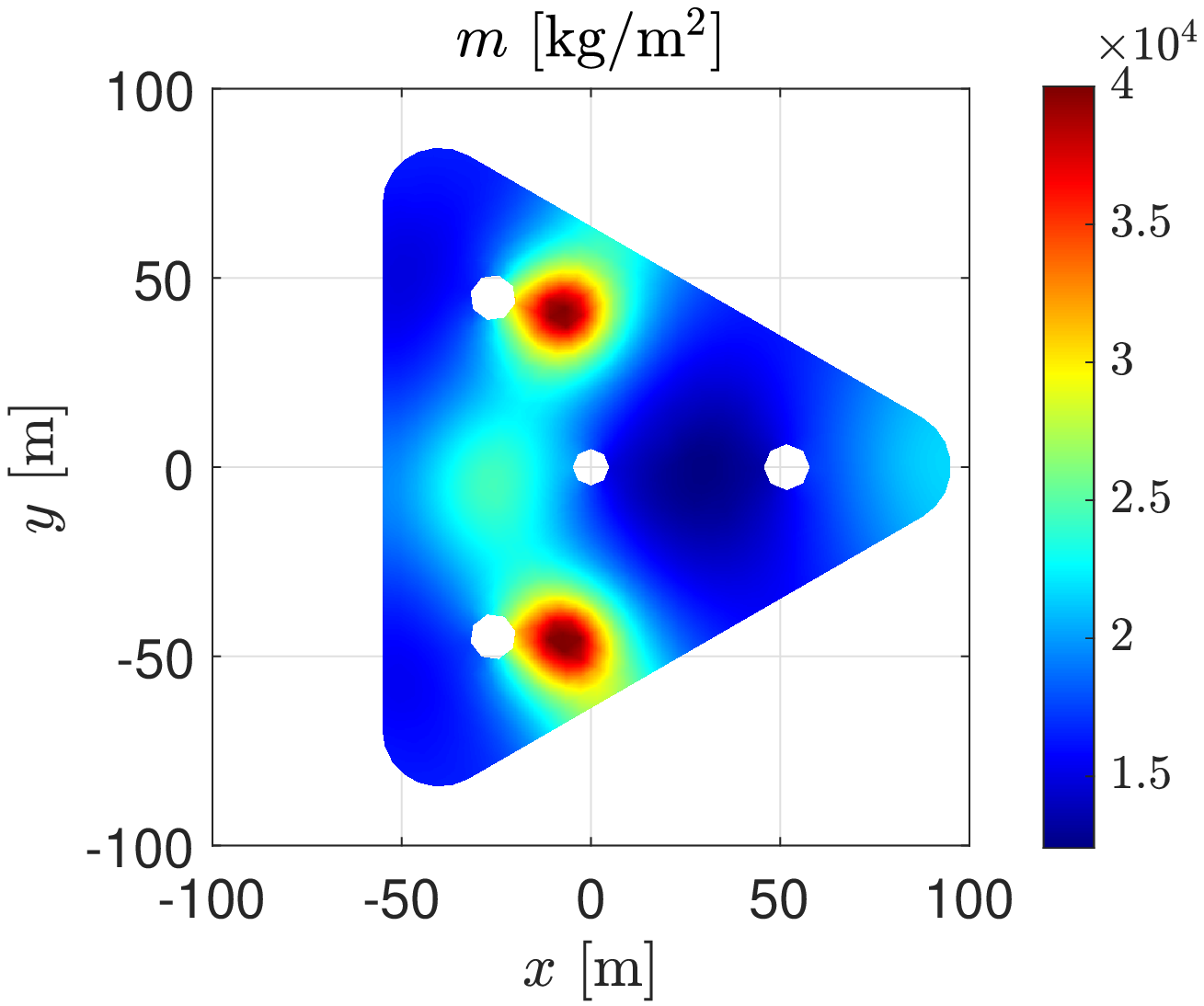}}
        &
        \subfloat[]{\label{fig:turb h}
        \includegraphics[width=0.32\textwidth,trim={15 0 30 -10}]{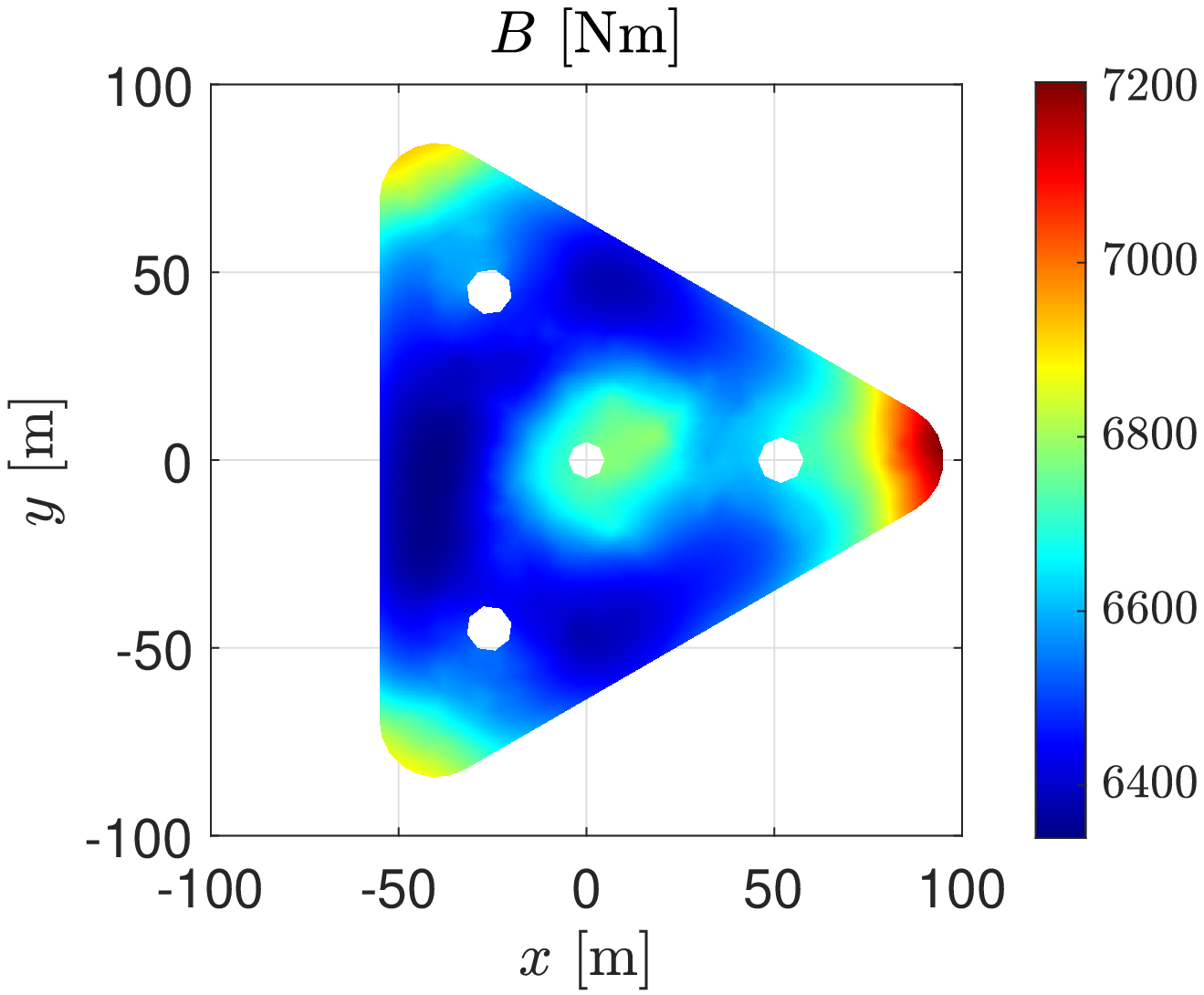}}
    \end{tabular}
    \caption{(a)~The shape of the floating wind turbine VolturnUS-S considered in the study. (b)~and~(c) show respectively amplitude and phase of the pressure applied; (d)~and~(e) surface mass and tension of the membrane; (f) costs over iterations for the cases of membrane and plate control; (g)~and~(h) surface mass and flexural stiffness of the plate.
    }
\end{figure}

In the following, we apply the three control strategies to the semi-submersible floating wind turbine depicted in figure~\ref{fig:turb a}, that is a standard defined by the European project~\cite{IEAWind}. The geometry and the mass parameters are computed according to the reference, so the stiffness and mass matrices $K$ and $M$ are defined by equation~\eqref{eq:M and K matrices}.
The non-null eigenvalues of $M^{-1}K$ approximate the resonances of the floating structures, that correspond to the time periods of $T_1=\SI{13.47}{\second}$  for heave motion, $T_2=\SI{36.33}{\second}$ for pitch and $T_3=\SI{36.30}{\second}$ for roll. By selecting the excitation period to be $\SI{13}{\second}$, the system is excited about resonance at $\omega=\SI{0.48}{\radian\per\second}$.
The seabed is set to be \SI{100}{\meter} dept and the regularisation
parameters are set as $\alpha_u=\alpha_v=\SI{1e-4}{}$, $\beta_u=\beta_v=\SI{1e-1}{}$ for both the membrane and plate controls.

Figures from~\ref{fig:turb b} to~\ref{fig:turb h} show the control actions applied and in table \ref{tab:bottle} one can appreciate the oscillation reduction in terms of amplitude, phase and minima of the cost functionals according to the three problems defined above.

\begin{table}
$ \begin{array}{ccccc}
\toprule
\{\bb X\} [\SI{}{\meter}/\SI{}{\degree}\angle\, \SI{}{\radian}] & \text{Uncontrolled} & \text{Pressure-controlled} & \text{Membrane-controlled} & \text{Plate-controlled} \\
\midrule
x &        0.613\,\angle-1.63  &  0.292\,\angle-1.63  &  0.056\,\angle -1.87  &  0.055\,\angle -1.38 \\
y &        0.000\,\angle-2.35  &  0.000\,\angle-1.53  &  0.053\,\angle\,1.61  &  0.012\,\angle\,1.59 \\
z &        0.708\,\angle-3.11  &  0.121\,\angle-3.10  &  0.070\,\angle -1.86  &  0.057\,\angle -0.36 \\
\theta_x & 0.000\,\angle-4.68  &  0.000\,\angle-3.81  &  0.221\,\angle -3.61  &  0.148\,\angle\,2.72 \\
\theta_y & 0.103\,\angle-4.27  &  0.015\,\angle-3.29  &  1.706\,\angle -2.96  &  0.484\,\angle\,3.39 \\
\midrule
J                                    & \SI{4.38e-1}{} & \SI{2.14e-1}{} & \SI{1.81e0}{} & \SI{1.82e0}{} \\
\frac{1}{2}\{\bb X\}^\dag C\{\bb X\} & \SI{4.38e-1}{} & \SI{5.00e-2}{} & \SI{5.86e-3}{} & \SI{3.22e-3}{}\\
\bottomrule
\end{array} $
\caption{Comparison between the turbine motion in case of the three different control strategies; the values describing rotation about $z$-axis is not reported because, it is almost zero in every case and derive from numerical errors only.}
\label{tab:bottle}
\end{table}



\section{Conclusions}

In this paper, we introduced a control-theoretic framework to modify the interaction between water waves and floating objects by acting on their surroundings through both active and passive control mechanisms. The active control consists of applied pressure which is optimally modulated in space to reduce the target's motion. The resulting OCP is linear-quadratic and it is solved efficiently with a {\it one-shot} approach. This active mechanism greatly reduces the oscillations but it may be difficult to realize in practice. As a consequence, we introduced two passive control mechanisms which are themselves floating objects with tunable properties, a membrane and a thin plate. By optimally modulating in space their physical properties, we showed how these mechanisms can still reduce the oscillations by two and one order of magnitudes, respectively. The control problem is bilinear in both cases and a system of first-order optimality conditions is derived using variational calculus. The resulting system is then solved iteratively.

On one hand, the results in this paper pave the way for a general and rigorous vibration control strategy which also provides effective design guidelines for reducing hydrodynamic excitation on arbitrary objects.
On the other hand, the framework we introduced is rather general and several other problems, such as cloaking, energy harvesting, etc., can be solved by properly modifying the cost functional.
Future work includes removing the small wave amplitude assumption and therefore considering higher order terms, resulting in a fully nonlinear problem where nonlinearities arise not only from the full Navier-Stokes equations but also from the fluid surface which is itself an unknown of the problem.

\vspace{15pt}


\bibliographystyle{RS}
\bibliography{main}

\end{document}